\pdfoutput=1

\documentclass[9pt, twocolumn, twoside]{IEEEtran}

\IEEEoverridecommandlockouts

\usepackage{amsthm}

\usepackage{amsmath,amssymb}
\usepackage{graphicx}
\usepackage{verbatim}
\usepackage{subfig}
\usepackage{subeqnarray}
\usepackage{setspace}
\usepackage{ifthen}
\usepackage{multicol}
\usepackage{float}
\usepackage{cite}
\usepackage[usenames, dvipsnames]{color}
\usepackage{caption}
\usepackage{comment} 
\usepackage{enumerate}
\usepackage{mathtools}
\usepackage{soul}
\usepackage{flushend}

\usepackage{color}

\DeclareMathOperator*{\argmax}{arg\,max}

\def\ba{\begin{array}}
\def\ea{\end{array}}
\newcommand{\beq}{\begin{equation}}
\newcommand{\eeq}{\end{equation}}
\newcommand{\bq}{\begin{eqnarray}}
\newcommand{\eq}{\end{eqnarray}}
\newcommand{\bqn}{\begin{eqnarray*}}
\newcommand{\eqn}{\end{eqnarray*}}
\newcommand{\bee}{\begin{enumerate}}
\newcommand{\eee}{\end{enumerate}}
\newcommand{\bi}{\begin{itemize}}
\newcommand{\ei}{\end{itemize}}
\newcommand{\ii}{\textbf{i}}

\newcommand{\gc}{\overline{g}_1}
\newcommand{\mingcf}{\min\left\lbrace\gc ,f_{12}\right\rbrace}
\newcommand{\mingf}{\min\left\lbrace\overline{g} ,f\right\rbrace}
\newcommand{\ac}{\alpha_\gamma}
\newcommand{\ad}{\alpha_\delta}

\newcommand{\TT}{\tau}
\newcommand{\Ps}{P^{s}}
\newcommand{\ps}{p^{s}}

\renewcommand\Re{\operatorname{\emph{Re}}}

\makeatletter
\newtheorem*{rep@theorem}{\rep@title}
\newcommand{\newreptheorem}[2]{%
\newenvironment{rep#1}[1]{%
 \def\rep@title{#2 \ref{##1}}%
 \begin{rep@theorem}}%
 {\end{rep@theorem}}}
\makeatother

\newtheorem{theorem}{Theorem}
\newreptheorem{theorem}{Theorem}
\newtheorem{lemma}{Lemma}
\newtheorem{proposition}{Proposition}
\newreptheorem{proposition}{Proposition}

\newboolean{showcomments}
\setboolean{showcomments}{true}

\newcommand{\bose}[1]{  \ifthenelse{\boolean{showcomments}}
{ \textcolor{blue}{(Bose says:  #1)}} {}  }
\newcommand{\chris}[1]{\ifthenelse{\boolean{showcomments}}
{ \textcolor{magenta}{(Chris says: #1)} } {} }
\newcommand{\babak}[1]{\ifthenelse{\boolean{showcomments}}
{ \textcolor{green}{(Babak says:  #1)}}{}}

\newcommand{\new}[1]{\ifthenelse{\boolean{showcomments}}
{ \textcolor{Red}{#1}}{}}

\graphicspath{{./Figures/}}

\usepackage[hidelinks]{hyperref}
\usepackage{comment}
\usepackage{thm-restate, enumerate }

\definecolor{darkred}{RGB}{150,0,0}
\definecolor{darkgreen}{RGB}{0,150,0}
\definecolor{darkblue}{RGB}{0,0,200}

\begin{document}
\title{ \bf Optimal Placement of Distributed Energy Storage in Power Networks}
\author{Christos Thrampoulidis, ~\IEEEmembership{Student Member,~IEEE,}
Subhonmesh Bose, ~\IEEEmembership{Student Member,~IEEE,} 
and Babak Hassibi ~\IEEEmembership{Fellow,~IEEE}.
\thanks{Emails: {\tt(cthrampo, boses, hassibi)@caltech.edu}. 
This work was supported in part by NSF under
grants CCF-0729203, CNS-0932428 and CCF-1018927, NetSE grant CNS 0911041, Office of
Naval Research MURI grant N00014-08-1-0747, Caltech
Lee Center for Adv. Net., ARPA-E grant DE-AR0000226, 
Southern Cali. Edison, Nat. Sc. Council of Taiwan, R.O.C. grant NSC 101-3113-P-008-001,
Resnick Inst., Okawa Foundation and Andreas Mentzelopoulos Scholarships for the Univ. of Patras.}}

\maketitle
\thispagestyle{empty}
\pagestyle{empty}

\begin{abstract}
We formulate the optimal placement, sizing and control of storage devices in a power network to minimize generation costs with the intent of load shifting. We assume deterministic demand, a linearized DC approximated power flow model and a fixed available storage budget. Our main result proves that when the generation costs are convex and nondecreasing, there always exists an optimal storage capacity allocation that places zero storage at generation-only buses that connect to the rest of the network via single links. This holds regardless of the demand profiles, generation capacities, line-flow limits and characteristics of the storage technologies. Through a counterexample, we illustrate that this result is not generally true for generation buses with multiple connections. For specific network topologies, we also characterize the dependence of the optimal generation cost on the available storage budget, generation capacities and flow constraints. 

\end{abstract}


 \section{Introduction}
\label{sec:intro}

\subsection{Motivation}

Energy storage technologies have been argued as \emph{``critical to achieving national energy policy objectives and creating a modern and secure electric grid system."} \cite{NAATBatt}. They have many potential applications in power networks, e.g., see \cite{chu2012, iecWhitePaper} for a detailed survey. At faster time scales (seconds to minutes), storage can be used to reduce variability of renewable sources of energy like wind or solar \cite{AEl11, xu2010load, budischak2012cost, sawindpower}. At slower time scales (over hours), it can be used for load shifting \cite{eyer2010energy, rastler2010electricity}, i.e., generate when it is cheaper and use storage dynamics to follow the demand. 
Though still very expensive, storage devices based on pumped hydro, compressed air, Lithium-ion based and other technologies have shown significant technical improvements and cost drops \cite{rastler2010electricity, Cal2020} over the last decade and are expected to play a central role in an efficient power system \cite{S04, nourai2007installation, denholm2010role, varaiya2011smart, eyer2010energy, Greenberger2011Online, lindley2010smart, NAATBatt, chu2012}.

Two natural questions to ask for storage are: (a) What is the optimal investment policy for storage? Where to place them, and how to size them? (b) Once installed, what is the optimal control policy for the storage as well as the generation schedule to minimize generation costs? In this paper, we formulate both problems for \emph{slower time-scales} in a common framework and present results on sizing such storage units in a network and a charging/ discharging policy for the installed units.

\subsection{Prior work}
Now, we provide a brief overview of the relevant literature. Optimal control policy for storage units has been extensively studied. While the authors in \cite{koutsopoulos2011optimal, su2011modeling, CLTX10} examine the control of a single storage device without a network, the authors in \cite{kanoria2011distributed, gayme_topcu2011} explicitly model the role of the networks in the operation of distributed storage resources. Storage resources at each node in the network are assumed to be known \emph{a priori} in these settings.

Sizing of storage devices has been studied by several authors, e.g., \cite{kraning2011, denholm09}  using purely economic arguments, without explicitly considering the network constraints of the physical system. Authors in \cite{su2011modeling,harshaoptimal} have looked at optimal sizing of  storage devices  in single-bus power systems for fast-time scales, 
 while Kanoria et al. \cite{ kanoria2011distributed} compute the effect of sizing of distributed storage resources on generation cost for specific networks.

The optimal storage placement problem on a general power network has been formulated and studied recently through simulations. The network imposes non-convex power-flow constraints that render such optimization problems NP-hard. These are handled through (a) linearization using DC approximation
\cite{purchala2005usefulness, stott2009}, or (b) a relaxation of the feasible sets using semidefinite programming \cite{bai2008, javad10, bose2012qcqp}. For the storage placement problem, Sj\"odin et al. in \cite{emma2011} uses the former, while Bose et al. in \cite{Bose2012storage} uses the latter.

\subsection{{Our Contribution}}

In this paper, we study the joint investment decision and control problem for storage devices in a power network. Our main contribution is the result in Theorem \ref{thm:main}: \emph{when minimizing a convex and nondecreasing generation cost with any fixed available storage budget over a slow time-scale of operation, there always exists an optimal storage allocation that assigns zero storage at nodes with only generation that connect via single transmission lines to the rest of the network. This holds for arbitrary demand profiles and other network parameters}.

First, we describe the salient features of our model. As in \cite{kanoria2011distributed,su2011modeling}, the investment decision problem is an infinite horizon problem. Hourly aggregate demands over large geographical locations often show periodicity \cite{caiso} and hence the optimization can be equivalently solved over one time period. The storage units are assumed to have finite capacities and ramp rates. Power exchanges with these devices suffer losses due to inefficiencies. Also, since we optimize the amount of storage placed on each bus, storage is assumed to be infinitely divisible. This is, however, not a limitation for our main result, as is explained in Section \ref{sec:discussion}. The generators have finite capacities with convex nondecreasing costs \cite{bergenBook,kanoria2011distributed,CLTX10}. The network has been modeled using linearized DC power-flow approximation \cite{Purchala2005} with finite line-flow capacities. This neglects reactive power over the network, defines voltage magnitudes to be at their nominal values at all buses and assumes the voltage phase angle differences between nodes to be small. Though this is a simplification of the full AC model of the power system with known limitations \cite{stott2009}, this approximation is widely used for analysis in optimal power flow \cite{stott2009,Purchala2005,Pandya08}, transmission expansion planning \cite{silva2006transmission} and electricity market operations \cite{rau2003issues,rau2003optimization,wuunifying}. The focus of this work is to derive structural properties of the storage placement problem using the linearized DC model as in \cite{emma2011, kanoria2011distributed}; this complements the studies without network models in \cite{CLTX10, su2011modeling} and simulation studies with a full AC model of the power flow equations \cite{gayme_topcu2011,Bose2012storage}. The result generalizes our work in \cite{CBBpes2013} and provides (partial) analytic justification of the observation made empirically in \cite{Bose2012storage, gayme_topcu2011, emma2011}: optimal storage allocation seldom places storage capacities at generator-only buses.

Next, we briefly discuss some of the qualifications of this work. First, to solve a complete storage investment strategy, we need a cost-benefit analysis of installing this new technology. In other words, the savings due to storage needs to be matched with the cost of installation and operation of such units on the grid. In this work, however, we only focus on minimizing cost of generation that estimates the potential benefits of storage.
Second, our main result applies to bulk storage on a slow time-scale and does not naturally generalize to scenarios with intermittent renewable generation. Dealing with fast time-scale variability of generation needs a stochastic control framework as in \cite{su2011modeling,kanoria2011distributed}; this, however, is not the focus of the current paper. 
Third, our main result characterizes storage allocation at generation-only buses that link to the rest of the network via a single transmission line. As shown in Section \ref{sec:multGen}, the result does not necessarily hold for generator-only buses with multiple links to the network. Also, it does not address the sizing or placement for any other kind of nodes in the network. We emphasize that this is a preliminary work on storage placement. To the best of our knowledge, this is the first theoretical result on this problem so far over a general power network. Our analysis suggests that there is a potential to exploit the rich underlying structure of this problem;  ongoing research aims at finding these properties and overcoming the limitations mentioned.

The paper is organized as follows. We formulate the optimal storage placement problem in Section \ref{problemFormulation}. The main result is stated and proven in Section \ref{sec:results}. A detailed discussion on the interpretation and some extensions of the result are presented in Section \ref{sec:discussion}.
We conclude the paper in Section \ref{conclusions}. 



\section{Problem formulation}
\label{problemFormulation}

Consider a power network that is defined by an undirected connected graph $\mathcal{G}$ on $n$ nodes (or buses) $\mathcal{N} = \{1, 2, \dots, n\}$. For two nodes $k$ and $l$ in $\mathcal{N}$, let $k \sim l$ denote that $k$ is connected to $l$ in $\mathcal{G}$ by a transmission line. We model time to be discrete and indexed by $t$. Now consider the following notation.
\begin{itemize}
	
\item $d_k(t)$ is the known real power demand at bus $k \in \mathcal{N}$ at time $t$. Hourly demand profiles often show diurnal variations \cite{pjm}, i.e., they exhibit cyclic behavior. Let $T$ time-steps denote the cycle length of the variation. In particular, for all $k \in \mathcal{N}$, $t \geq 0$, assume
$$
d_k (t + T) = d_k (t).
$$

\item $g_k(t)$ is the real power generation at bus $k \in \mathcal{N}$ at time $t$ and it satisfies
\bq
\label{eq:limits for g}
0 \leq g_k(t) \leq \overline{g}_k,
\eq
where, $\overline{g}_k$ is the generation capacity at bus $k$.\footnote{We do not include ramp constraints on the generators.} 


\item $c_k \left(g_k \right)$ denotes the cost of generating power $g_k$ at bus $k\in \mathcal{N}$. The cost of generation is assumed to be independent of time $t$ and depends only on the generation technology at bus $k$. Also, suppose that the function $c_k : \mathbb{R}^+ \to \mathbb{R}^+$ is non-decreasing
 and convex. These assumptions apply to commonly used cost functions in the literature \cite{bergenBook, javad10, gayme_topcu2011, bose2012qcqp}, e.g., convex and nondecreasing piecewise linear or quadratic ones.

\item The power $p_{kl}$ sent from bus $k$ towards bus $l$ for two nodes $k \sim l$ in $\mathcal{G}$ is limited by thermal and stability constraints as
\beq
| p_{kl} (t) | \leq {f}_{kl},
\label{eq:flowLim}
\eeq
where $f_{kl}$ is the capacity of the corresponding line.

\begin{figure}[h]
\centering
	\includegraphics[width=0.30\textwidth]{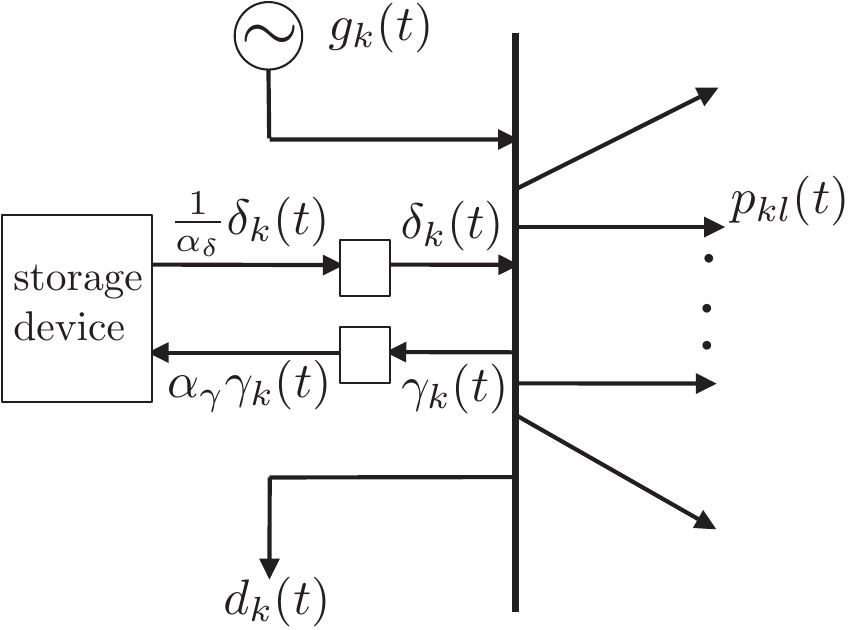}
	\caption{Power balance at node $k\in\mathcal{N}$.}
		\label{fig:powBalance}	
\end{figure}

\item $\gamma_k(t)$ and $\delta_k(t)$ are the average charging and discharging powers of the storage unit at bus $k \in \mathcal{N}$ at time $t$, respectively. The energy transacted over a time-step is converted to power units by dividing it by the length of the time-step. This transformation conveniently allows us to formulate the problem in units of power \cite{Bose2012storage}. Let $0 < \ac, \ad \leq 1$ denote the charging and discharging efficiencies, respectively of the storage technology used, i.e., the power flowing in and out of the storage device at node $k\in\mathcal{N}$ at time $t$ is $\ac\gamma_k(t)$ and $\frac{1}{\ad}\delta_k(t)$, respectively \cite{su2011modeling, korpaas2003operation}. The roundtrip efficiency of this storage technology is $\alpha=\ac\ad \leq 1$.

\item $s_k(t)$ denotes the storage level at node $k \in \mathcal{N}$ at time $t$ and $s_k^0$ is the storage level at node $k$ at time $t=0$. From the definitions above, we have
\begin{align}\label{eq:s_k(t)}
s_k(t) = s_k^0 + \displaystyle\sum_{\tau=1}^{t}{\left( \ac\gamma_k(\tau) - \frac{1}{\ad}\delta_k(\tau) \right)}.
\end{align}
For each $k \in \mathcal{N}$, assume $s_k^0 = 0$, so that the storage units are empty at installation time.

\item $b_k \geq 0$ is the storage capacity at bus $k$. Thus, $s_k (t)$ for all $t$ satisfies
\beq
\label{eq:storage profile constraints}
0 \leq s_k(t) \leq b_k.
\eeq

\item $h$ is the available storage budget and denotes the total amount of storage capacity that can be installed in the network. Our optimization algorithm decides the allocation of storage capacity $b_k$ at each node $k\in\mathcal{N}$ and thus, we have\footnote{Note that we do not restrict the capacity sizes $b_k, k\in\mathcal{N}$ \emph{a priori}. The problem formulation can be extended to include any linear constraints on $b_k$'s.}
\beq
\label{eq:sum(bk)<h}
\sum_{k \in \mathcal{N}} b_k \leq h.
\eeq

\item The charging and discharging rates of storage device at node $k \in \mathcal{N}$ are bounded above by ramp limits; these limits are assumed to be proportional to the installed storage capacity at node $k$, i.e.,\footnote{{Note that $\epsilon_\gamma$ and $\epsilon_\delta$ are specific to a storage technology. We consider the installment of one kind of storage over the network. Though we present our results with one storage technology, it can be generalized to the joint placement of installing multiple storage types with individual storage budgets for each technology.}}
\begin{subequations}\label{eq:gammaDelta}
\begin{align}
0 \leq \gamma_k(t) \leq \epsilon_\gamma b_k, \\
0 \leq \delta_k(t) \leq \epsilon_\delta b_k,
\end{align}
\end{subequations}
where $\epsilon_\gamma\in(0, \frac{1}{\ac} ]$ and $\epsilon_\delta\in(0,\ad]$ are fixed constants.
\end{itemize}

Balancing power that flows in and out of bus $k\in \mathcal{N}$ at time $t$, as shown in Figure \ref{fig:powBalance}, we have:
\beq 
\label{eq:power balance}
g_k(t)-d_k(t)-\gamma_k(t) + \delta_k(t) = \displaystyle\sum_{l \sim k} p_{kl}(t).
\eeq 

The power flow $p_{kl}$ from bus $k$ to bus $l$ relates to the voltages at the buses through Kirchoff's law. Since power is quadratic in voltage, the power flow equations introduce quadratic equalities that render most optimization problems over power networks nonconvex and hence hard to solve and analyze. The role of nonconvexity in power flow optimization has been widely studied in the power system literature 
; see \cite{o2012iv,frank2012optimal} surveys. Nonconvex optimizations, in general, are hard to solve and difficult to analyze. To make the model amenable to analysis, one option is to linearize the power flow equations around an operating point. Such a linearization technique popularly used in the literature is the DC approximation \cite[Ch.\ 9]{grainger1994power} \cite[Ch.\ 6]{andersson2004modelling} \cite{ purchala2005usefulness}; for completeness, we discuss it in Appendix \ref{app:DCApp}. 
 In this model, the transmission losses (resistances in transmission lines) and reactive power flows in the network are ignored. Specifically, suppose $B_{kl}$ is the susceptance of the transmission line joining buses $k$ and $l$ and $\theta_k(t)$ is the voltage phase angle at bus $k \in \mathcal{N}$ at time $t$. Then, using DC approximation, it can be shown that 
\beq
\label{eq:line flow eqn}
p_{kl} (t) = B_{kl}\left[\theta_{k}(t)-\theta_{l}(t)\right].
\eeq
Though we ignore all transmission losses for presenting our result, we generalize it in Section \ref{sec:losses} to include losses. The loss model used is discussed further in Appendix \ref{app:DCApp}.


Optimally placing storage over an infinite horizon is  equivalent to solving this problem over a singe cycle, provided the state of the storage levels at the end of a cycle is the same as its initial condition \cite{Bose2012storage}. Thus, for each $k \in \mathcal{N} $, we have
\beq
\label{eq:sum(r)=0}
\sum_{t=1}^{T}{ \left( \ac\gamma_k(t) - \frac{1}{\ad}\delta_k(t) \right) } = 0.
\eeq
For convenience, denote $[T]:=\left\lbrace 1,2,\dots,T \right\rbrace$. Using the above notation, we define the following optimization problem.\\
\noindent
\textbf{Storage placement problem $P$:}
\begin{equation*}\label{SP}
\begin{aligned}
\underset{}{\text{minimize}}  & \ \ \ \sum_{k \in \mathcal{N}} \sum_{t=1}^{T}{ c_k \left(g_k(t)\right) } \\
\text{over} & \ \ \ (g_k(t),\gamma_k(t),\delta_k(t),p_{kl}(t),b_k, \theta_k(t)), 
\\ & \ \   \ k\in\mathcal{N},  \  k\sim l, \ t\in [T],
\\
\text{subject to} & \ \ \ \eqref{eq:limits for g} - \eqref{eq:sum(r)=0},\\
\end{aligned}
\end{equation*}
where, \eqref{eq:limits for g} represents generation constraints, \eqref{eq:flowLim}, \eqref{eq:power balance} represent power flow constraints, \eqref{eq:storage profile constraints},\eqref{eq:sum(bk)<h},\eqref{eq:gammaDelta},\eqref{eq:sum(r)=0} represent the constraints imposed on the charging/discharging control policy of the energy storage devices, and \eqref{eq:line flow eqn} represents the DC approximated Kirchoff's laws. For the power network, $\theta_k (t), k \in \mathcal{N}$ and $p_{kl}(t), k \sim l \text{ in } \mathcal{G}$ are state variables, while $g_k(t),\gamma_k(t),\delta_k(t), b_k$ are controllable inputs to the system.

Given the demand profiles and network parameters, $P$ can be efficiently solved to define the optimal investment decision strategy for sizing storage units at different buses, the economic dispatch of the various generators and the optimal control policy of the installed storage units.
%

Now, restrict attention to network topologies where each bus either has generation or load but not both\footnote{An intermediate bus (one that has no generation or load) is modeled as a load bus. Any losses associated with the node is included as a load.}. Partition the set of buses $\mathcal{N}$ into two groups $\mathcal{N}_G$ and $\mathcal{N}_D$ where they represent the generation-only and load-only buses respectively and assume $\mathcal{N}_G$ and $\mathcal{N}_D$ are non-empty. For any subset $\mathcal{K}$ of $\mathcal{N}_G$, define the following optimization problem.

\noindent
\textbf{Restricted storage placement problem $\Pi^{\mathcal{K}}$:}
\begin{equation*}\label{AP}
\begin{aligned}
\underset{}{\text{minimize}}  & \ \ \ \sum_{k \in \mathcal{N}} \sum_{t=1}^{T}{ c_k \left(g_k(t)\right) } \\
\text{over} & \ \ \ (g_k(t),\gamma_k(t),\delta_k(t),\theta_k(t),p_{kl}(t),b_k), 
\\ & \ \   \ k\in\mathcal{N},  \  k\sim l, \ t\in [T],
\\
\text{subject to} & \ \ \ \eqref{eq:limits for g} - \eqref{eq:sum(r)=0},\\
& \ \ \ b_k = 0, \quad k\in\mathcal{K}.
\end{aligned}
\end{equation*}
Problem $\Pi^{\mathcal{K}}$ corresponds to placing no storage at the (generation) buses of the network in subset $\mathcal{K}$. We study the relation between the problems $P$ and $\Pi^{\mathcal{K}}$ in the rest of the paper. 

\begin{figure}[h!]
\centering
	\includegraphics[width=0.28\textwidth]{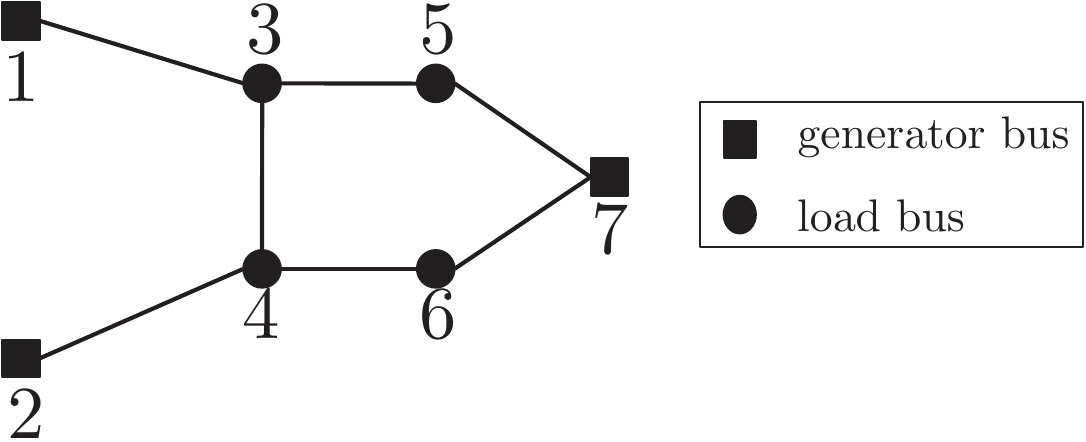}
	\caption{A sample network.}
		\label{fig:sample network}	
\end{figure}

We say bus $k \in \mathcal{N}$ has \emph{a single connection} if it has exactly one neighboring node $l \sim k$. Similarly, a bus $k \in \mathcal{N}$ has \emph{multiple connections} if it has more than one neighboring node in $\mathcal{G}$. We illustrate the notation using the network in Figure \ref{fig:sample network}. $\mathcal{N}_G = \{ 1, 2, 7 \}$ and $\mathcal{N}_D = \{ 3, 4, 5, 6 \}$. Buses 1 and 2 have single connections and all other buses in the network have multiple connections.  

\section{Main Result}
\label{sec:results}

For a subset $\mathcal{K} \subseteq \mathcal{N}_G$, let $p_*$ and ${\pi_*}^\mathcal{K}$ be the optimal values for problems $P$ and $\Pi^\mathcal{K}$, respectively. Now, we are ready to present the main result of this paper.

\begin{theorem}
\label{thm:main}
Let $\mathcal{K}\subseteq \mathcal{N}_G$ be a subset of generation nodes that have single connections.
Consider the storage placement problem $P$, the restricted storage placement problem $\Pi^{\mathcal{K}}$ and their respective optimal costs $p_*$ and $\pi_*^\mathcal{K}$. 
 If $P$ is feasible, then $\Pi^\mathcal{K}$ is feasible and $p_* = \pi_*^\mathcal{K}$.
\end{theorem}
Problem $P$, in general, may have multiple optimal solutions, but Theorem \ref{thm:main} proves that there  
\emph{always} exists an optimal allocation of storage capacities that places \emph{no} storage at any subset of generation buses with single connections, regardless of the demand profiles, generation capacities, line-flow limits and characteristics of the storage technologies.


Next, we make the following remarks about the result.
(a) Notice that we have restricted our attention to generator buses in $\mathcal{K}$ that have single connections only. This does not generalize to the case if $\mathcal{K}$ includes generator buses with multiple connections; see Section \ref{sec:multGen} for an example. (b) Storage capacity allocation at each bus has been assumed to be infinitely divisible, i.e., each $b_k, k\in \mathcal{N}$ that satisfies the budget constraint $\sum_{k \in \mathcal{N}} b_k \leq h$ in \eqref{eq:sum(bk)<h} is feasible. But it might be impractical to implement an optimal allocation with arbitrarily small storage capacities. This, however, is not a limitation for the result in Theorem \ref{thm:main} as it only specifies zero storage capacities at some buses and does \emph{not} characterize storage sizes at others. (c) In our formulation, we assume perfect knowledge of the entire demand profile. The result in Theorem \ref{thm:main}, however, holds true for \emph{any} demand profile as long as the storage placement problem $P$ is feasible.

Before presenting the proof, we provide some intuition behind the result. Consider a generator bus $k$ that has a single connection to node $l$ in the network. First, we solve the storage placement problem $P$ for this network. Suppose this results in some storage capacity installed at bus $k$ with some charging/ discharging profile $\gamma_k^* (t), \delta_k^*(t), \ t \in [T]$. To construct a solution of the restricted storage problem $\Pi^{\{k\}}$, a natural idea to explore is to \emph{shift} this storage from bus $k$ to bus $l$ and operate it with the optimal control policy ($\gamma_k^* (t), \delta_k^*(t), \ t \in [T]$) obtained from the solution of $P$. This shift can be done, provided that the optimal generation profile $g_k^*(t), t \in [T]$, itself, defines a feasible flow over the  transmission line, i.e., $g_k^* (t) \leq f_{kl}$ for all $t \in [T]$. The key insight to prove this fact is that at the time instant where $g_k^*(t)$ is at its maximum, the storage at bus $k$ cannot be charging. If it was indeed charging, one could generate less and charge less at the same time. 
In what follows, we formalize this argument. 

\subsection{Proof of Theorem \ref{thm:main}}
We only prove the case where the round-trip efficiency is $\alpha < 1$, but the result holds for $\alpha = 1$ as well. Assume $P$ is feasible throughout. For any variable $z$ in problem $P$, let $z^*$ be the value of the corresponding variable at the optimum. In our proof, we use the following technical result.
\begin{lemma}\label{lemma:meanVal}
Suppose $\phi:\mathbb{R} \rightarrow \mathbb{R}$ is convex. Then, for any $ x_1 < x_2$ and $0\leq \eta \leq (x_2-x_1)$:
\begin{align*}
\phi( x_1  +\eta ) + \phi( x_2  - \eta ) \leq \phi( x_1 ) +  \phi( x_2 ).
\end{align*} 
\end{lemma}
\begin{IEEEproof}
Applying Jensen's inequality to the convex function $\phi(\cdot)$, we have
\begin{align*}
\left( 1 - \frac{\eta}{x_2 - x_1} \right) \phi(x_1) + \left( \frac{\eta}{x_2 - x_1}\right) \phi(x_2) \ &\geq \ \phi(x_1 + \eta),\\
\left( \frac{\eta}{x_2 - x_1} \right) \phi(x_1) + \left( 1 - \frac{\eta}{x_2 - x_1}\right) \phi(x_2) \ &\geq \ \phi(x_2 - \eta).
\end{align*}
The result follows from adding the inequalities above.
\end{IEEEproof}

Consider node $k \in \mathcal{K}$ and $k \sim l$. Node $l$ is uniquely defined as $k$ has a single connection. Problem $P$, in general, has multiple optima. In the following result, we characterize only a subset of these optima.
\begin{lemma}
\label{lemma:gf}
There exists an optimal solution of $P$ such that for all $t \in [T]$ and all $k \in \mathcal{K},  l\sim k$, 
\begin{enumerate}[(a)]
\item \label{lemma.1} $g^*_k(t)\gamma^*_k(t)\delta^*_k(t) = 0$,
\item \label{lemma.2} $g_k^*(t) \leq f_{kl}$.
\end{enumerate}
\end{lemma}
%

\begin{IEEEproof}The feasible set of problem $P$ is a bounded\footnote{ 
Without loss of generality, let bus $1$ be the slack bus and hence $\theta_1(t)=0$ for all $t\in[T]$. Boundedness of the set of feasible solutions of $P$ then follows from the relations in  \eqref{eq:limits for g},  \eqref{eq:flowLim}, \eqref{eq:sum(bk)<h}, \eqref{eq:gammaDelta} and \eqref{eq:line flow eqn}.
} polytope and the objective function is a continuous convex function. Hence the set of optima of $P$ is a convex compact set \cite{boyd_book}. Now, with every point in the set of optimal solutions of $P$, consider the function $\sum_{k \in \mathcal{K}, t\in [T]} \left( \gamma_k (t) + \delta_k (t) \right)$. This is a linear continuous function on the compact set of optima of $P$ and hence attains a minimum. Consider the optimum of $P$ where this minimum is attained. We prove parts (a) and (b) in Lemma \ref{lemma:gf} for this optimum.  

(a) Suppose, on the contrary, we have $g^*_k(t_0)>0$, $\gamma^*_k(t_0)>0$ and $\delta^*_k(t_0)>0$ for some $t_0 \in [T]$. Define 
\begin{align*}
\Delta g' := \min\left\{ (1-\alpha)\gamma_k^*(t_0) \ , \frac{1-\alpha}{\alpha}\delta_k^*(t_0) \ , g_k^*(t_0) \right\}.
\end{align*}
Note that $\Delta g'>0$. Now, for bus $k$, construct modified generation, charging and discharging profiles $\tilde{g}_k(t), \tilde{\delta}_k(t), \tilde{\gamma}_{k}(t), t \in [T]$ that differ from ${g}^*_k(t), {\delta}^*_k(t), {\gamma}^*_{k}(t)$ only at $t_0$ as follows:
\begin{align*}
\tilde{g}_k(t_0) & := g^*_k(t_0) - \Delta g ', 
\end{align*}
\begin{align*}
\tilde{\gamma}_k(t_0) & := \gamma^*_k(t_0)-\frac{1}{1-\alpha}\Delta g ', \\
\tilde{\delta}_k(t_0)  & := \delta_k^*(t_0) - \frac{\alpha}{1-\alpha}\Delta g '.
\end{align*}
Note that, for all $t\in[T]$, the storage level $s_k(t)$ and the power $p_{kl}(t)$ flowing from bus $k$ to bus $l$ remain unchanged throughout. It can be checked that the modified profiles define a feasible point of $P$. Since $c_k(\cdot)$ is non-decreasing, we have $c_k\left( \tilde{g}_k(t_0) \right) \leq c_k\left( g^*_k(t_0) \right)$ and hence the additivity of the objective in $P$ over $k$ and $t$ implies that this feasible point has an objective function value of at most $p_*$. It follows that this feasible point defines an optimal point of $P$. Also, we have  $\tilde{\gamma}_k(t_0) +  \tilde{\delta}_{k}(t_0)  < {\gamma}^*_k(t_0) +  {\delta}^*_{k}(t_0)$ and thus, this optimum of $P$ has a strictly lower $\sum_{k \in \mathcal{K}, t\in [T]} \left( \gamma_k (t) + \delta_k (t) \right)$, contradicting our hypothesis. This completes the proof of Lemma \ref{lemma:gf}(\ref{lemma.1}).


(b) If $g^*_k(t)=0$ for all $t\in [T]$, then $g_k^*(t) \leq f_{kl}$ trivially holds. Henceforth, assume $\max_{t\in [T]}{g^*_k(t)}> 0$, and consider $t^+\in [T]$, such that $g_k^*(t^+) = \max_{t\in[T]}g^*(t)$.


If $\gamma_k^*(t^+) = 0$, then,
\begin{align}\label{eq:stop}
\max_{t\in[T]}g_k^*(t) &= g_k^*(t^+) \notag 
= \underbrace{p^*_{kl}(t^+)}_{ \leq f_{kl} } + \underbrace{\gamma^*_k(t^+)}_{ =0 } - \underbrace{\delta^*_k(t^+)}_{ \geq 0 } \notag 
 \leq f_{kl},
\end{align}
and Lemma \ref{lemma:gf}(\ref{lemma.2}) holds.


Now, suppose that $\gamma_k^*(t^+)>0$ and hence $\delta_k^*(t^+)=0$ from Lemma \ref{lemma:gf}(\ref{lemma.1}). Since storage charges by an amount $\alpha_\gamma \gamma_k^*(t^+) > 0$, we have $s_k^*(t^+) > 0$. Also, $s_k^*(T) = s_k^0 = 0$ by hypothesis and hence the storage at node $k$ discharges from $s_k^*(t^+)$ to zero in $[t^+ + 1, T]$. Let $t^-$ be the first time instant after $t^+$ when the storage device at bus $k$ discharges, i.e.
\begin{equation}\label{eq:t1}
t^- := \min{ \left\{ t \in [t^++1,T] \ \Big| \ \ac\gamma_k^*(t) - \frac{1}{\ad}\delta^*_k(t) < 0 \right\}  }.
\end{equation}
Thus, $\delta_k^*(t^-)>0$. Now, we argue that $g_k^*(t^-) = g_k^*(t^+)$. Since $g_k^*(t^+) = \max_{t \in [T]} g_k^*(t)$, clearly $g_k^*(t^-) \leq g_k^*(t^+)$. Suppose this inequality is strict.
Then we show how to construct an optimum of $P$ with a lower 
$\sum_{k \in \mathcal{K}, t\in [T]} \left( \gamma_k (t) + \delta_k (t) \right)$ to contradict our hypothesis.
Define
 \begin{align*}
 \Delta g := \min \Big\{ \gamma^*_k(t^+),~ \frac{1}{\alpha}\delta_k^*(t^-),~  g_k^*(t^+),~ \frac{1}{\alpha}\left(g_k^*(t^+) - g_k^*(t^-)\right)\Big\}. 
 \end{align*}
Observe that $\Delta g > 0$. Construct the modified generation, charging and discharging profiles at node $k$,  $\tilde{g}_k(t), \tilde{\delta}_k(t), \tilde{\gamma}_{k}(t)$ that differ from ${g}_k^*(t), {\delta}_k^*(t), {\gamma}_{k}^*(t)$ only at $t^+$ and $t^-$ as follows:
$$
\tilde{g}_k(t^+) := g^*_k(t^+) - \Delta g, \quad \tilde{g}_k(t^-) := g^*_k(t^-) + \alpha \Delta g, 
$$
$$
\tilde{\gamma}_k(t^+) := \gamma^*_k(t^+) - \Delta g, \quad \tilde{\gamma}_k(t^-) := \gamma^*_k(t^-) ,
$$
$$
\tilde{\delta}_k(t^+) := \delta^*_k(t^+) = 0, \quad \tilde{\delta}_k(t^-) := \delta_k^*(t^-) - \alpha \Delta g.
$$
Also, define the modified storage level $\tilde{s}_k(t)$ using $\tilde{\gamma}_k(t)$ and $\tilde{\delta}_k(t)$. To provide intuition to the above modification, we essentially generate and store less at time $t^+$ by an amount $\Delta g$. This means at a future time $t^-$, we can discharge $\alpha \Delta g$ less from the storage device and hence have to generate $\alpha \Delta g$ more to compensate. From the definition of $\Delta g$, it follows that for $t= t^+, t^-$, we have $0 \leq \tilde{g}_k(t) \leq \overline{g}_k$, $0 \leq \tilde{\gamma}_k(t) \leq \epsilon_\gamma b^*_k$, and $0 \leq \tilde{\delta}_k(t) \leq \epsilon_\delta b^*_k$. Also, the line flows $p_{kl}(t)$ remain unchanged. For the storage levels, it can be checked that
\begin{align*}
0 \leq s^*_k(t^+-1) \leq \tilde{s}_k(t) \leq s_k^*(t) \leq b^*_k, & \ \text{for } t \in [t^+, t^--1], \\
\tilde{s}_k(t) = s_k^*(t), & \  \text{otherwise}.
\end{align*}
This proves that the modified profiles define a feasible point for $P$. Also, we have
\begin{subequations}\label{eq:compCost}
\begin{align}
 &c_k\left( \tilde{g}_k(t^+)  \right) + c_k\left( \tilde{g}_k(t^-)  \right) \notag\\
\label{eq:compCost.1}  & \qquad \leq c_k\left( g^*_k(t^+) - \alpha \Delta g  \right) + c_k\left( g^*_k(t^-) + \alpha \Delta g   \right) \\
\label{eq:compCost.2} & \qquad \leq c_k\left( {g}^*_k(t^+)  \right) + c_k\left( {g}^*_k(t^-)  \right),
\end{align}
\end{subequations}
where \eqref{eq:compCost.1} follows from the non-decreasing nature of $c_k(\cdot)$ and \eqref{eq:compCost.2} follows from our assumption $g_k^*(t^-) < g_k^*(t^+)$ and Lemma \ref{lemma:meanVal}. 
The modified profiles $\tilde{g}_k(t), \tilde{\delta}_k(t), \tilde{\gamma}_{k}(t)$ are feasible in $P$ with an objective value at most $p_*$. Thus, they define an optimum of $P$. Also, 
\begin{align*}
&\tilde{\gamma}_k(t^+) + \tilde{\gamma}_k(t^-) + \tilde{\delta}_k(t^+) + \tilde{\delta}_k(t^-)\\
& \quad = {\gamma}_k^*(t^+) + {\gamma}_k^*(t^-) + {\delta}_k^*(t^+) + {\delta}_k^*(t^-) - \underbrace{(1 + \alpha) \Delta g}_{>0}.
\end{align*}
This implies that this optimum of $P$ has a lower $\sum_{k \in \mathcal{K}, t\in [T]} \left( \gamma_k (t) + \delta_k (t) \right)$, contradicting our hypothesis. Hence, we have $g_k^*(t^-) = g_k^*(t^+) = \max_{t \in [T]}g_k^*(t)$.

Finally,
$$\max_{t\in[T]}g_k^*(t) = g_k^*(t^-) \notag 
= \underbrace{p^*_{kl}(t^-)}_{ \leq f_{kl} } + \underbrace{\gamma^*_k(t^-)}_{ =0 } - \underbrace{\delta^*_k(t^+)}_{ > 0 }  
 \leq f_{kl}.
$$
This completes the proof of the lemma.
\end{IEEEproof}

To prove Theorem \ref{thm:main}, consider the optimal solution of $P$ that satisfies Lemma \ref{lemma:gf}(\ref{lemma.2}). For $k \in \mathcal{K}$, $g_k^*(t)$ itself defines a feasible flow over the line joining buses $k$ and $l$, where $l$ is the unique neighboring node of $k$ in graph $\mathcal{G}$. Now, transfer the storage device at bus $k$ to bus $l$. In particular, define a new storage capacity $\hat{b}_l := b_k^* + b_l^*$ and operate it with a charging (and discharging) profile $\hat{\gamma}_l (t) = \gamma_k^*(t) + \gamma_l^* (t)$ (and similarly for $\hat{\delta}_l (t)$). For node $k$, define the new voltage phase angle $\hat{\theta}_k(t) := \theta^*_k(t)+\dfrac{1}{B_{kl}} ( \gamma_k^*(t) - \delta_k^*(t) )$. The power flow from bus $k$ to bus $l$ is then given as $\hat{p}_{kl}(t) := g_k^*(t)$. These profiles define a feasible point of $\Pi^{ \{ k \}}$ with an objective value of $p_*$. Combining with the fact that $p_* \leq \pi^{ \{ k \}}_*$, we conclude $p_* = \pi^{ \{ k \}}_*$. Finally, we do this successively for each $k \in \mathcal{K}$ to obtain $p_* = \pi^{  \mathcal{K} }_*$.



\section{Discussion and extensions}
\label{sec:discussion}

Here, we discuss our main result in more detail. We begin by considering cases where Theorem \ref{thm:main} is helpful to the network planner in Section \ref{sec:application}. Then in Section \ref{sec:losses}, we present an extension of our result to the case with losses in the network. We comment on the importance of convexity in the problem formulation in Section \ref{sec:concaveCost} and finally explore storage placement at buses with multiple connections in Section \ref{sec:multGen}.

\subsection{Applicability of Theorem \ref{thm:main}}
\label{sec:application}

\begin{figure}[H]

\centering
\subfloat[]{
	\centering
\hspace{-0pt}	\vspace{-5pt}
	\includegraphics[width=0.115\textwidth]{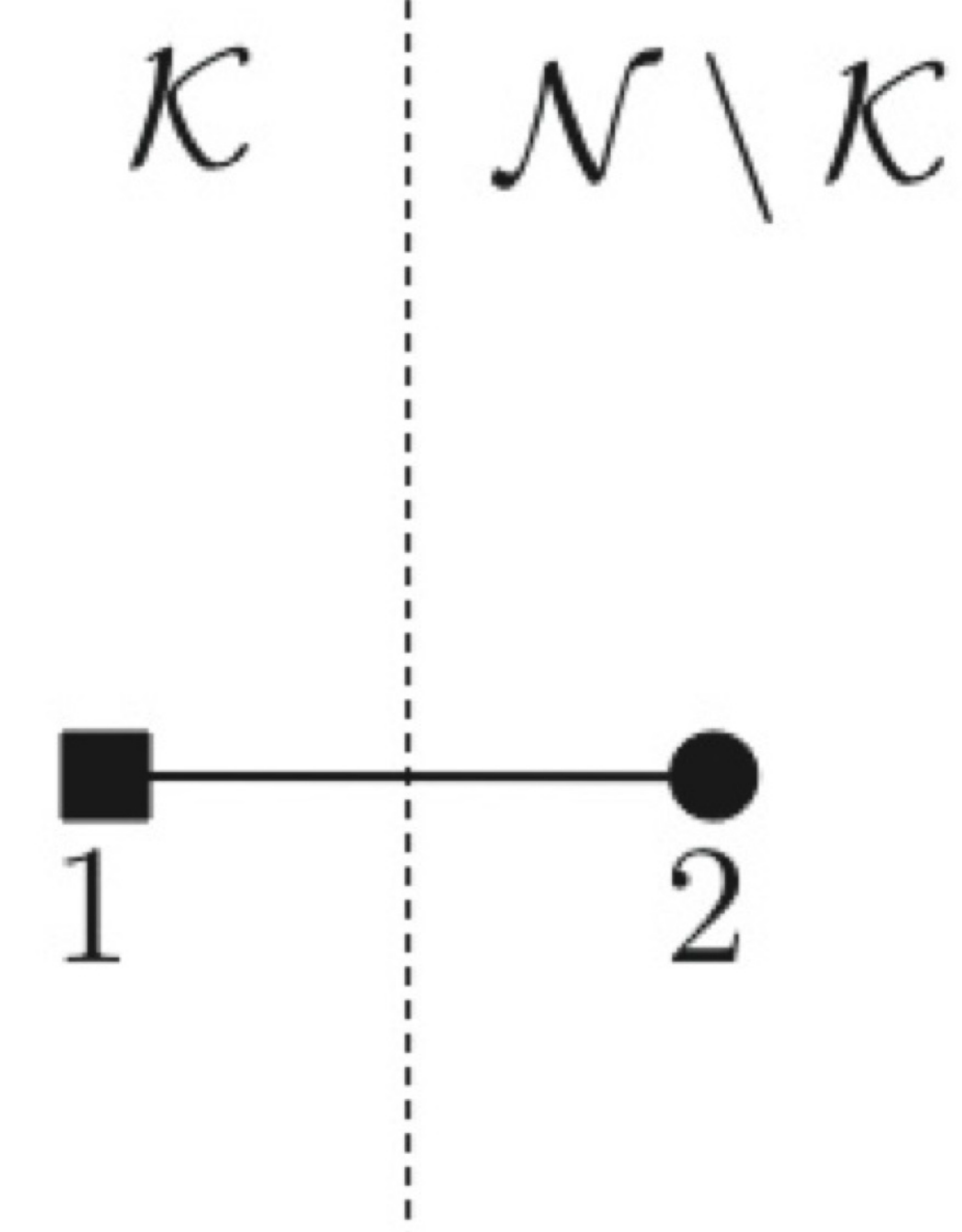}
	\label{fig:ex1}
}
\subfloat[]{
\vspace{2pt}
	\centering
\hspace{15pt}	\includegraphics[width=0.195\textwidth]{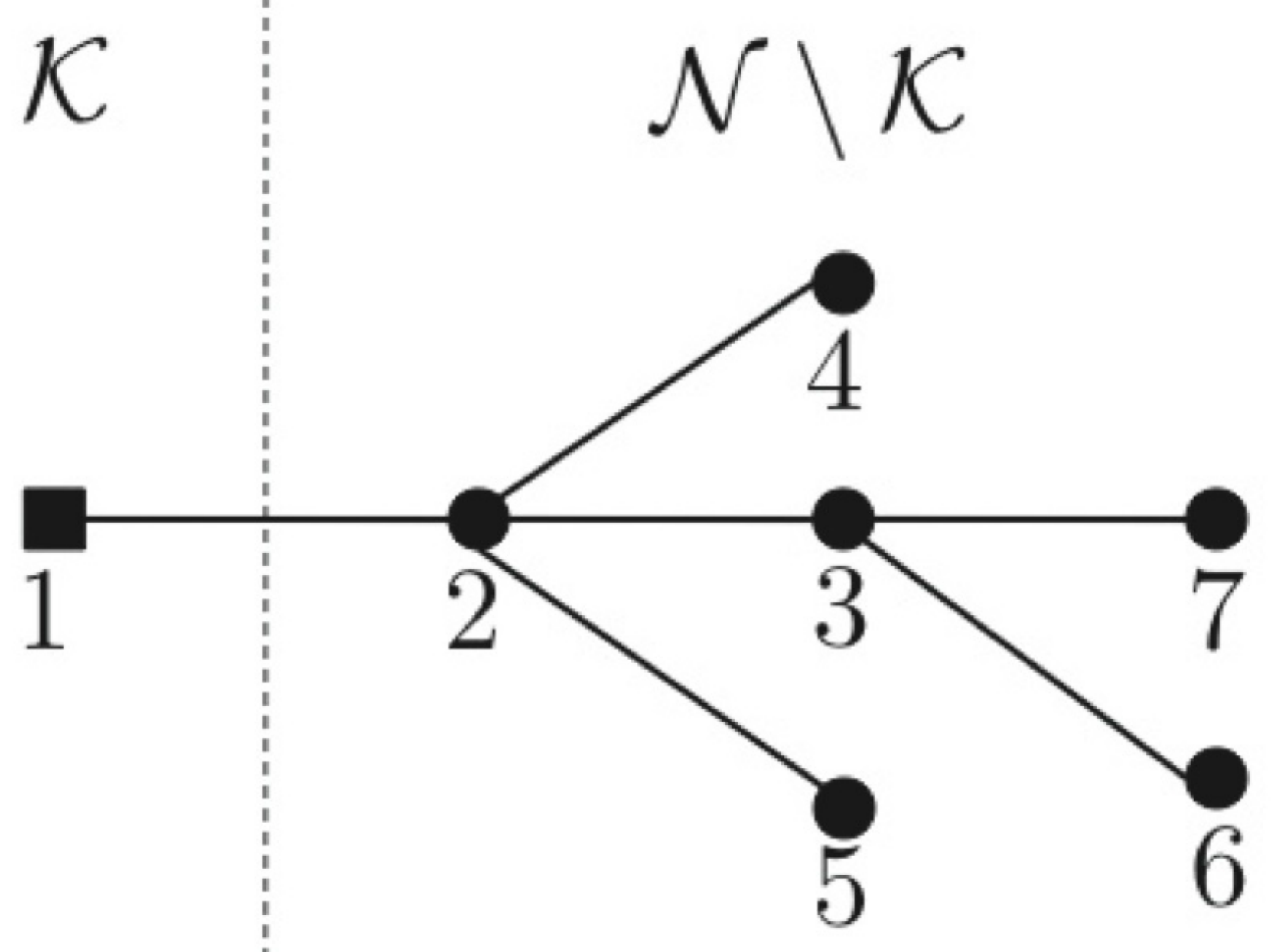}
	\label{fig:ex2}
}	
\caption{Examples of power networks (a) Single generator single load system (b) A radial network.}
\label{fig:examples}

\end{figure}

Figure \ref{fig:examples} depicts a few power networks, where Theorem \ref{thm:main} applies, i.e., network topologies with generator buses that have single connections. While this may seem quite restrictive, in practice, many networks have generators of this type. The single generator single load case in Figure \ref{fig:ex1} models topologies where generators and loads are geographically separated and are connected by a transmission line, e.g., see \cite{6281676}. This is common where the resources for the generation technology (like coal or natural gas) are available far away from where the loads are located in a network. Figure \ref{fig:ex2} is an example of a radial network, i.e., an acyclic graph. Most distribution networks conform to this topology\footnote{Two assumptions in our model hold for transmission networks but not strictly for distribution networks: (a) Resistances in distribution lines are not negligible and hence DC approximation does not generally apply \cite{stott2009}, (b) Three different types of loads, namely, constant power, constant current and constant impedance loads show different behavior in distribution networks \cite{bergenBook}; but in aggregate, demands can be modeled as constant power loads in transmission networks, as in IEEE distribution feeders \cite{IEEEdistr}.}, e.g., see \cite{IEEEdistr, bose2012qcqp}. Also, isolated transmission networks, e.g., the power network in Catalina island \cite{xu2010load} are radial in nature. 

Next, we discuss how Theorem \ref{thm:main} defines an investment strategy that is robust to many system parameters. Our result suggests that it remains optimal not to place any storage at buses in set $\mathcal{K}$ even if the demand profiles, generation capacities, line flow capacities or admittances in the network change. Consider the example in Figure \ref{fig:ex1}. Suppose the line flow capacity is larger than the peak value of the demand profile, i.e., $f_{12} \geq \max_{t \in [T]} d_2(t)$. It can be checked that placing all the available storage at the generator bus is also an optimal solution. If at a later time during the operation of the network, the demand increases such that the peak demand surpasses the line capacity, this placement of storage no longer remains optimal and requires new infrastructure for storage to be built on the demand side to avoid load shedding. If, however, we use the optimum as suggested by the problem $\Pi^\mathcal{K}$ and place all storage on the demand side from the beginning, then this placement not only can accommodate the change in the demand, but it also remains optimal under the available storage budget. To explore another such direction, suppose another generator is built to supply the load in Figure \ref{fig:ex1}. From Theorem \ref{thm:main}, it follows that we still do not need storage allocation at bus 1 even with the extended network. 

\subsection{Modeling losses in the network}
\label{sec:losses}

\begin{figure}[h!]
\centering
	\includegraphics[width=0.32\textwidth]{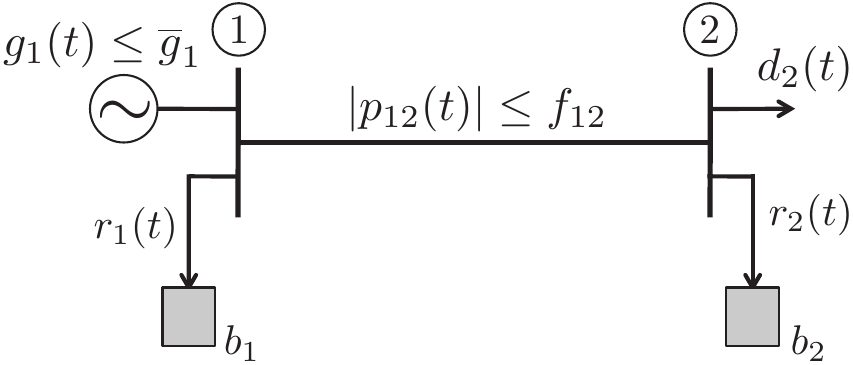}
	\caption{Single generator single load network. Available storage budget is $h \geq b_1+b_2$.}
	\label{fig:SGSL}
\end{figure}

Our problem formulation uses linearized DC approximation for the power flow equations. This approximation ignores all transmission losses in the network. Here, we explore one popular way to incorporate losses and generalize the result in Theorem \ref{thm:main}. To simplify the presentation, consider the single generator single load network shown in Figure \ref{fig:SGSL}. Generator at bus 1 is connected to a load at bus 2 using a single line, i.e., $\mathcal{K} = \mathcal{N}_G = \{1\}$ and  $\mathcal{N}_D = \{2\}$. We make further simplifications: assume $\ac = \ad = 1$, i.e., the storage devices are perfectly efficient. Hence, define for $t\in[T]$ and $k = 1,2$,
$$r_k(t) := \gamma_k (t) - \delta_k (t).$$

Recall that $p_{12}$ is the power injected at bus 1 towards bus 2. However, it suffers a loss before reaching bus 2. As detailed in Appendix \ref{app:DCApp}, losses in such a network can be approximated to be quadratic in the power sent, i.e., loss $ \approx \xi p_{12}^2$, where $\xi > 0$ is some positive constant that depends on the impedance of the transmission line. Thus, power received at bus 2 is $p_{12} - \xi p_{12}^2$. Then balancing power on both buses, we get for $t \in [T]$, 
\begin{align}
\label{eq:quadLoss}
p_{12}(t) = g_1(t) - r_1(t), \ \ \ p_{12}(t) - \xi p_{12}^2(t) = d_2 (t) + r_2(t).
\end{align}

Let the storage placement problem with losses incorporated be ${P}^L$. Following the definition of  ${\Pi}^{\{1\}}$, define the restricted storage placement problem with losses as ${\Pi}^{\{1\},L}$; this is essentially the problem ${P}^L$ with the extra constraint $b_1 = 0$. Let the optimal costs of these problems be ${p}_*^L$ and ${\pi}^{\{1\},L}_*$, respectively. With this notation, we have the following result.
\begin{proposition}
\label{prop:loss}
The storage placement problems with losses ${P}^L$ and ${\Pi}^{\{1\},L}$ satisfy ${p}_*^L = {\pi}^{\{1\},L}_*$.
\end{proposition}
In what follows, we provide a proof sketch; details are deferred to Appendix \ref{sec:lossesProof}.
Perhaps the first thing one notices about the problems ${P}^L$ and ${\Pi}^{\{1\},L}$ is that  they are nonconvex due to the quadratic equality in \eqref{eq:quadLoss}. Modify the problems ${P}^L$ and  ${\Pi}^{\{1\},L}$ to their \emph{convex relaxations}, where the second equality in  \eqref{eq:quadLoss} is changed to $p_{12}(t) - \xi p_{12}^2(t) \geq d_2 (t) + r_2(t)$. Call these relaxations as problems $\hat{P}^L$ and  $\hat{\Pi}^{\{1\},L}$, respectively. Let their optimal costs be $\hat{p}_*^L$ and $\hat{\pi}^{\{1\},L}_*$, respectively. Using the set inclusion relations among the feasible sets of the programs  ${\Pi}^{\{1\},L}$, ${P}^L$ and $\hat P^{L}$, it is easy to argue that
\begin{align}
\label{eq:lossProof}
\hat{p}_*^L  \leq  {p}_*^L  \leq  {\pi}^{\{1\},L}_* .
\end{align}
Then the proof proceeds in two steps. 
First, we show $\hat{p}_*^L =  \hat{\pi}^{\{1\},L}_*$ and then prove that the relaxation $\hat{\Pi}^{\{1\},L}$ is tight, i.e. $\hat{\pi}^{\{1\},L}_* = {\pi}^{\{1\},L}_*$. 
Using \eqref{eq:lossProof}, it is straightforward to see that these two statements imply Proposition \ref{prop:loss}. 




We briefly discuss our result in Proposition \ref{prop:loss} here. First, notice that the optimization problem $P^L$ is non-convex due to the constraint in \eqref{eq:quadLoss}. Consequently, it is computationally hard to solve. However, $\hat\Pi^{\{1\},L}$ considers a convex relaxation of this nonconvex constraint and $\hat\Pi^{\{1\},L}$ is a convex program that often admits an efficient solution. In addition, we have $p_*^L = \hat{\pi}^{\{1\},L}_* $ and thus $\hat\Pi^{\{1\},L}$ provides a computationally tractable way of \emph{exactly} solving $P^L$. Second, notice that Proposition \ref{prop:loss} considers losses that are quadratic in the power sent. The proof technique generalizes to the case with any \emph{convex} function of the power sent. Third, we have presented Proposition \ref{prop:loss} for a two-node network only. This idea can be generalized to a network with losses to derive a result similar to Theorem \ref{thm:main}, i.e., in a network with quadratic (or convex) losses, there exists an optimal storage placement with zero storage at generators with single connections.



\subsection{Effect of concave cost functions}
\label{sec:concaveCost}
Theorem \ref{thm:main} assumes a nondecreasing \emph{convex} cost of generation; this is commonly found in practice, e.g., the costs of coal-based generators are often increasing quadratic functions \cite{bergenBook}. Convexity is sufficient for our result to hold, but may not be necessary. However, the following example shows that the theorem need not generalize for arbitrary non-decreasing functions.
Consider the network in Figure \ref{fig:SGSL} and
let the cost of generation at bus 1 be a concave cost function 
$c_1(g_1)= 2 g_1$, if $0\leq g_1\leq 5$ and $c_1(g_1)=g_1 + 5$ otherwise. With $T=2$, let the load bus have a demand profile  $d_2 = ( 5, 5)$ and $f_{12}=5$ connecting them. Further let $h=1$, $\alpha=1$, $\epsilon_\gamma=\epsilon_\delta=1$ and $\overline{g}_1=8$. All quantities are in per units. It can be checked that the optimal generation profile of $\Pi^{\{1\}}$ is $(5,5)$, thus, $\pi_*^{\{1\}}=20$. On the other hand, the generation profile $(6,4)$ is feasible for $P$. Hence, $p_*\leq 19 < \pi_*^{\{1\}}$. We also remark that when $c(\cdot)$ is not convex, $P$ and $\Pi^{\mathcal{K}}$ are not convex programs and, hence, cannot be solved efficiently.

\subsection{On generators with multiple connections}
\label{sec:multGen}
\begin{figure}[h!]

\centering
	\includegraphics[width=0.10\textwidth]{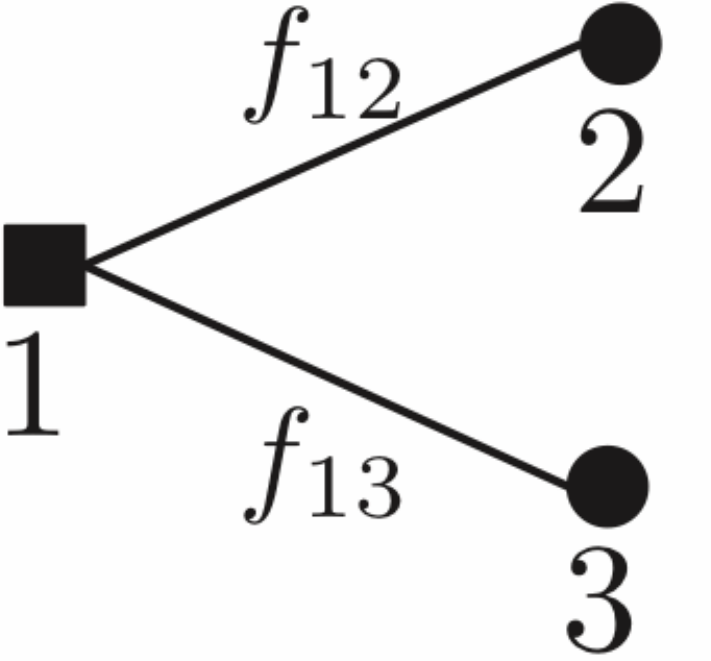}
	\caption{A 3-node network with a generator bus with multiple connections.}
	\label{fig:counter}
\end{figure}

Our result in Theorem \ref{thm:main} considers generator buses that have single connections only. A natural direction to generalize the result is to include generator buses with multiple connections. However, we show through an example that generator buses with multiple connections may not always have zero storage capacity in the optimal allocation. Consider a 3-node network as shown in Figure \ref{fig:counter}. All quantities are in per units. Let the cost of generation at node 1 be $c_1(g_1) = g_1^2$. Let $T=4$ and the demand profiles at nodes 2 and 3 be $d_2=(9,10,0,10)$ and $d_3=(0,10,9,10)$.
Also, suppose that the line and generation capacities are $f_{12} = f_{13} = 9.5$ 
and the available storage budget is $h=5$. Finally, assume no losses and ignore the ramp constraints in the charging and discharging processes, i.e. $\alpha=1$ and $\epsilon_\gamma=\epsilon_\delta=1$. The optimal storage allocation $(b_1^*, b_2^*, b_3^*)$ for the two problems $P$ and $\Pi^{\{1\}}$ is $(4, 0.5, 0.5)$ and $(0, 2.5, 2.5)$, respectively. Also, the optimal generation profile $g_1^*(t),t=1,2,3,4$ for the two problems can be computed to be $(14, 15, 14, 15)$ and $(12, 17,12, 17)$, respectively. Thus, $p_* = 842 < \pi_*^{\{1\}} = 866$.

%
%

We provide some intuition behind the design of the counterexample. First, notice that if demands at buses 2 and 3 are multiples of each other, i.e., $d_2(t) = \phi d_3(t)$ for some constant $\phi \geq 0$ (and so are the line capacities), the 3-node network can be roughly thought of as two single-generator-single-load systems with nodes $(1,2)$ and $(1,3)$, respectively and one can again prove that the generator node does not need any storage capacity in an optimal allocation. We formalize this statement in Proposition \ref{prop:proportionalStar} in Appendix \ref{app:partial}. Thus to expect $b_1^* \neq 0$ in $P$, we consider demand profiles that show opposite trends. Second, if $h=\infty$ or $h$ is the minimum value required for a feasible flow, we show in Proposition \ref{prop:hmin_starN} in Appendix \ref{app:partial} that there exists an optimal point with $b_1^* = 0$. Hence, we consider a storage budget that is in the middle. Third, note that if line capacities are large, then an optimal allocation with $b_1^*=0$ trivially exists. Thus, we construct $f_{12}=f_{13} =9.5$ for which $P$ and $\Pi^{\{1\}}$ are feasible but the network is congested. This illustrates some key directions to look at for characterizing cases where $b_1^*=0$ for generator buses with multiple connections.

\section{Concluding remarks}
\label{conclusions}

In this paper, we formulate the optimal storage placement problem for load shifting at slow time scales. We show in Theorem \ref{thm:main} that generator nodes with single connections get zero storage capacity at optimal allocation; this holds regardless of the demand profiles, generation capacities, line-flow limits and characteristics of the storage technologies. 
The counterexample in Section \ref{sec:multGen} shows that such a general result does not hold beyond the settings in Theorem \ref{thm:main}. However, notice that the demand profiles of the two demand nodes considered in this counterexample show opposite
 variations. When demands vary similarly, e.g., in Appendix \ref{app:partial} (specifically Proposition \ref{prop:proportionalStar}), it can be shown that one does not need storage installed at the generator node for the same network. This suggests that one would expect to discover structural results when demand profiles and/ or network topologies are restricted to a certain class. Characterizing these classes is a natural direction to explore.

In our formulation, we have neglected any fixed installment cost of storage; we have only minimized the total generation cost over every cycle of operation. It is easy to argue that the result extends to the case where such installment costs are linear in the storage capacity installed. The role of general (possibly nonlinear) installment costs on investment decisions would be an interesting are of study.

As a final remark, we have only focussed on minimizing the cost of conventional generation. The primary use of storage in this work is load shifting, i.e., arbitrage between power consumed at different times through storage rather than following the load variations with generation. However, an important application of storage is to mitigate intermittency of stochastic renewable generation at faster time-scales. We would like to pursue the implications and extensions of this paper to such scenarios.




\section{Acknowledgements}
The authors gratefully acknowledge Prof. K. Mani Chandy (Caltech), Mr. Paul DeMartini (Resnick Institute, Caltech), Christophe Hennekinne (Cambridge) and all the reviewers for their comments.

\bibliography{Powerbib}

\appendices

\section{DC Approximation and losses}
\label{app:DCApp}
In this section, we derive the linearized DC approximation and the loss model for the power flow equations from Kirchoff's laws. We start with introducing some notation. Define $\ii := \sqrt{-1}$ and for any complex number $z$, let $z^H$ define its conjugate.

Recall that the power network is defined by an undirected connected graph $\mathcal{G}$ on $n$ nodes (or buses) $\mathcal{N} = \{1, 2, \dots, n\}$. Let $V_j$ be the voltage at bus $j$. Since, the power network works in sinusoidal alternating current mode, the voltages and currents are represented as complex numbers\footnote{This is often referred to as the phasor representation \cite{bergenBook} that essentially represents the time-varying signal in the Fourier domain.}. For convenience, let $V_j := v_j e^{\ii \theta_j}$, where $ v_j \geq 0$ is the magnitude and $\theta_j$ is the argument of $V_j$. The transmission line joining buses $j$ and $k$ has an admittance $y_{jk} = G_{jk} - \ii B_{jk}$. Usually, $B_{jk}$ is positive since most transmission lines are inductive in nature. Also, let there be shunt elements $y_{jj} = G_{jj} - \ii B_{jj}$ associated with each bus. For the shunt element, however, $B_{jj}$ is usually non-negative since shunt elements are generally capacitive\footnote{We refer the reader to \cite{bergenBook} for the one phase equivalent lumped circuit model of power systems.}. 

From Kirchoff's laws, it then follows that the apparent power injection (generation - demand) $p_j + \ii q_j$ at bus $j$ satisfies 
\begin{align}
& p_j  + \ii q_j  \notag\\
& = V_j {{\left(V_j y_{jj}\right)^H} + \sum_{k\sim j} V_j {\left( (V_j - V_k) y_{jk} \right)^H}} \notag\\
& = v_j^2 \left(y_{jj}^H + \sum_{k \sim j} y_{jk}^H \right) -  \sum_{k\sim j} v_j v_k e^{\ii (\theta_j - \theta_k)} y_{jk}^H \notag
%
%
%
\end{align}
\begin{align}
& = \left [ v_j^2 G_{jj}  +  v_j^2  \sum_{k\sim j} G_{jk} - \sum_{k\sim j} v_j v_k G_{jk} \cos (\theta_j - \theta_k) \right.\notag\\
& \left. \qquad \qquad +  \sum_{k\sim j} v_j v_k B_{jk} \sin (\theta_j - \theta_k) \right] \notag\\ 
& \quad + \ii \left [ - v_j^2 B_{jj} + v_j^2  \sum_{k\sim j} B_{jk}  -  \sum_{k\sim j} v_j v_k G_{jk} \sin(\theta_j - \theta_k) \right. \notag\\
& \left. \qquad \qquad - \sum_{k\sim j} v_j v_k B_{jk} \cos(\theta_j - \theta_k) \right]. \label{eq:pjqj}
\end{align}

\subsection{DC Approximation}
Now, we introduce the DC approximation \cite[Ch.\ 9]{grainger1994power}, \cite[Ch.\ 6]{andersson2004modelling}, \cite{bergenBook, purchala2005usefulness, stott2009dc}:
\begin{enumerate}
\item Resistances in transmission lines are small compared to the inductances and hence ignored, i.e., $G_{jk} \approx 0$. Shunt reactances are similarly ignored, i.e., $B_{jj} \approx 0$.
\item Voltage magnitudes are maintained close to their nominal values. Measured in per units, $ v_j \approx 1$ for all $j \in \mathcal{N}$.
\item Voltage phase angle differences across a transmission line are small, i.e., for $j \sim k$, we have $\sin(\theta_j - \theta_k) \approx (\theta_j - \theta_k)$ and $\cos(\theta_j - \theta_k) \approx 1$.
\end{enumerate}

Using the above approximation in \eqref{eq:pjqj}, we get
\begin{align}
p_j  + \ii q_j  = G_{jj}  +  \sum_{k\sim j} B_{jk} (\theta_j - \theta_k). 
\label{eq:pjqj.DC}
\end{align}
Notice that the right hand side of \eqref{eq:pjqj.DC} is real and hence $q_j = 0$; in DC approximation, there is no reactive power flow. The loss associated with bus $j$ (or more accurately in its shunt element) is $G_{jj}$ and is usually included in the net real power demand at bus $j$.

\subsection{Modeling Losses}\label{sec:app_losses}
Now, we turn to modeling losses. We refer the reader to \cite{coffrin2012approximating, stott2009dc} for detailed discussions on incorporating losses in linearized DC approximated power flow equations. For any complex number $z$, let $| z |$ denote its magnitude and $\Re z$ denote the real part of $z$. For nodes $j \sim k$ in graph $\mathcal{G}$, notice that the current flowing from bus $j$ to bus $k$ is given by $I_{jk} := (V_j - V_k) y_{jk}$. 
%
The resistance on this line is $\Re (1/y_{jk}) = G_{jk}/(G_{jk}^2 + B_{jk}^2)$. Then the loss on this line is given by 
\begin{align*}
| I_{jk} |^2 \Re (1/y_{jk})  
& = | V_j - V_k |^2 |y_{jk}|^2 G_{jk}/(G_{jk}^2 + B_{jk}^2) \\
& = \left[ v_j^2 + v_k^2 - 2 v_j v_k \cos (\theta_j - \theta_k) \right]  G_{jk}.
\end{align*}
Using the approximations $v_j, v_k \approx 1$ and $\cos (\theta_j - \theta_k)\approx 1 - (\theta_j - \theta_k)^2/2$, we have 
\begin{align*}
| I_{jk} |^2 \Re (1/y_{jk})  
& \approx (\theta_j - \theta_k)^2  G_{jk}.
\end{align*}
Given the DC approximation, the power that flows from bus $j$ to $k$ is $p_{jk} = B_{jk} (\theta_j - \theta_k)$. Thus the loss incurred on the transmission line is given by 
$$(\theta_j - \theta_k)^2  G_{jk} = \left(\frac{G_{jk}}{B_{jk}^2}\right)p_{jk}^2 .$$

\section{Proof of Proposition \ref{prop:loss}}\label{sec:lossesProof}

Recall the definitions of problems $\hat P^L$ and $\hat \Pi^{\{1\},L}$ defined for the single generator single load network in Figure \ref{fig:SGSL}. As outlined in Section \ref{sec:losses}, the proof consists of two steps: (a) $\hat p^L_* = \hat \pi^{\{1\},L}_*$ and (b) $\hat \pi^{\{1\},L}_* = \pi^{\{1\},L}_*. $

\subsection{Proof of $\hat p^L_* = \hat \pi^{\{1\},L}_*$}\label{sec:hatp=hatpi}
From the size of the feasible sets, it is easy to check that $\hat p^L_* \leq \hat \pi^{\{1\},L}_*.$ We will show that equality holds in the above relation as follows: we start with an optimal solution of $\hat P^L$ and appropriately modify it to yield a  solution with the same generation profile (thus, the same cost) but such that it satisfies $r_1(t)=0$ (thus, feasible for $\hat \Pi^{\{1\},L}$). Denote the generation profile, storage profiles and storage capacities of such an optimal solution of $\hat P^L$ as $g_1^*(t)$, $r^*_1(t)$, $r^*_2(t)$, $b_1^*$ and $b_2^*$, respectively. Also, let $s^*_1(t)=\sum_{\tau=1}^{t} r_1^*(\tau)$ denote the optimal storage level at the generation node at time $t$. 
From the feasibility constraints of $\hat P^L$:
\begin{align}\label{eq:weKnow}
g_1^*(t) - r^*_1(t) - \xi (g^*_1(t) - r^*_1(t))^2 \geq d_2 (t) + r^*_2(t).
\end{align}
The main idea of the proof is similar to that of Theorem \ref{thm:main}: we ``shift" the storage with capacity $b_1^*$ from node $1$ to node $2$, i.e., the modified storage capacity at node 2 becomes $b_2^*+b_1^*$ and its associated storage profile becomes $r_2^*(t) + \tilde r(t), t\in[T]$. Notice that the added term $\tilde r(t)$ must be a feasible storage profile, i.e., it satisfies \eqref{eq:storage profile constraints} and \eqref{eq:sum(r)=0}; they are repeated here for convenience. For $t \in [T]$,
\begin{align}\label{eq:suf1}
0\leq \sum_{\tau=1}^{t} \tilde{r} (\tau)\leq b_1^* \quad \text{ and } \quad
\sum_{\tau=1}^{T}\tilde r(\tau) = 0.
\end{align}
Also, it must satisfy the power balance equation 
 \begin{align}\label{eq:needs} 
g_1^*(t) - \xi (g^*_1(t) )^2 \geq d_2 (t) + r^*_2(t) + \tilde r(t).
\end{align}
In the absence of losses ($\xi = 0$) in Theorem \ref{thm:main}, it was sufficient to choose $\tilde r(t) = r_1^*(t), t\in[T]$; this choice, however, need not satisfy \eqref{eq:needs} for $\xi>0$. In what follows, we show how to construct a suitable $\tilde r(t),t\in [T]$.

\begin{figure}[h!]
     \centering
\includegraphics[scale=0.38]{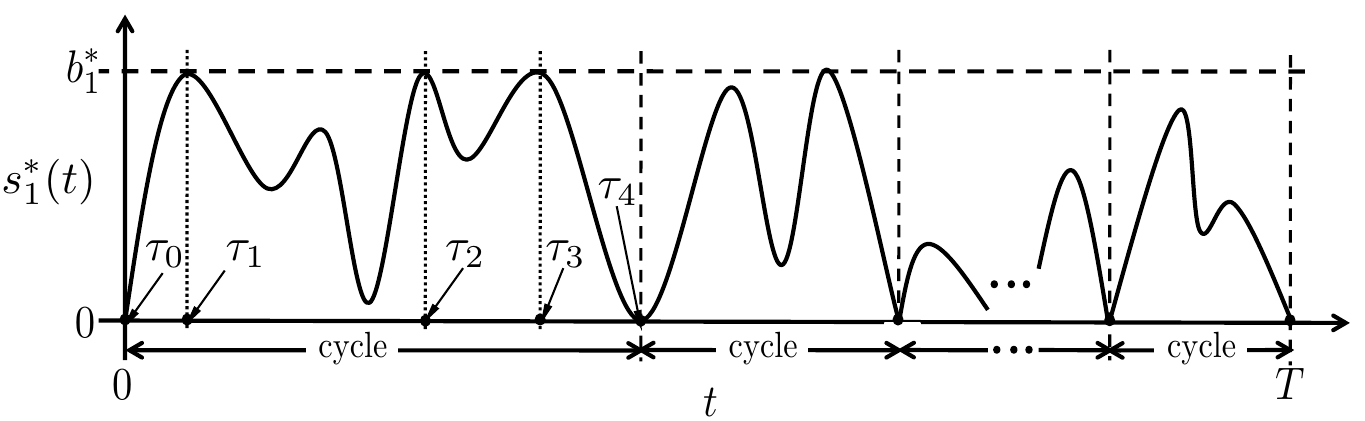}	
		\caption{Dividing the time horizon $[0,T]$ into ``cycles".}
		\label{fig:cycle}
\end{figure}

Consider a representative storage profile $r_1^*(t)$ in Figure \ref{fig:cycle}. We divide the time horizon $[0,T]$ into ``\emph{cycles}" which are time intervals with the property that $s^*_1(t) = 0$ is only at the start and end of the interval and strictly positive in between. 
We construct $\tilde r(t)$ that has the same set of cycles as $r_1^*(t)$, i.e.
\begin{align}\label{eq:stateCyc}
\sum_{t\in\text{ cycle }} \tilde r(t) = 0.
\end{align}
 It is easy to check that it suffices to construct $\tilde r(t)$ over one such cycle. To reflect this change, the relations in \eqref{eq:suf1} are modified as follows: the summations run over $\tau$ from the start of a cycle and $T$ is replaced by the end of the cycle. 

For convenience, define $[t_1, t_2] := \{ t_1,t_1+1,\ldots,t_2 \}$. Without loss of generality, we construct $\tilde{r}(t)$ over the first cycle $[\TT_0:=0,\TT_{m+1}]$ as in Figure \ref{fig:cycle}. Also, let $\TT_1 < \TT_2 < \ldots < \TT_m \in [\TT_0+1, \TT_{m+1}]$ be the sequence of time instants at which the storage at node $1$ is at its capacity\footnote{If $s_1^*(t) \neq b_1^*$ within the cycle, then $m=0$.}. Formally, we have
\begin{align*}
\begin{cases}
s_1^*(t) = 0 & \text{ for } t=\TT_0,\TT_{m+1}\\
s_1^*(t) = b^*_1& \text{ for } t=\TT_1,\ldots, \TT_{m}\\
0<s_1^*(t)<b_1^* & \text{ otherwise} .
\end{cases}
\end{align*}

Consider the optimal solution of $\hat P^L$ which attains the minimum of the function $\sum_{t\in[T]}\left( |r_1(t)|+g_1(t) \right)$ over the set of optimal solutions of $\hat P^L$.  We begin by showing a useful property of this optimal solution in the following result.
\begin{lemma}\label{lem:12}
Consider $t^+,t^-\in[\TT_0+1,\TT_{m+1}]$ such that $r^*_1(t^+)>0$ and $r^*_1(t^-) < 0$. 
\begin{enumerate}[(a)]
\item If $t^+<t^-$, then $g_1^*(t^+)\leq g_1^*(t^-)$.
\item If there exists $i=0,1,\ldots,m$ such that $\TT_i+1\leq t^-<t^+\leq \TT_{i+1}$, then
$g_1^*(t^+)\leq g_1^*(t^-)$.
\end{enumerate}
\end{lemma}
\begin{IEEEproof}
The proof mimics the proof of Lemma \ref{lemma:gf}; we omit some details for brevity.

(a) Suppose to the contrary $g_1^*(t^+)>g_1^*(t^-)\geq 0$.
Define 
\begin{align*}
\Delta g: = \min\Big\{ r^*_1(t^+), -r^*_1(t^-), &g_1^*(t^+)-g_1^*(t^-), \\~& \min_{t\in[t^++1,t^--1]} s^*_1(t)\Big\}.
\end{align*}
Then it can be checked that $\Delta g>0$. Construct a modified generation, charging and discharging profile at node $1$,  $\tilde{g}_1(t), \tilde{r}_1(t)$, that differ from ${g}^*_1(t), r^*_1(t)$ only at $t^+$ and $t^-$ as follows:
\begin{align*}
\tilde{g}_1(t^+) = g^*_1(t^+) - \Delta g, \quad \tilde{g}_1(t^-) = g^*_1(t^-) +  \Delta g, \\
\tilde{r}_1(t^+) = r^*_1(t^+) - \Delta g, \quad \tilde{r}_1(t^-) = r^*_1(t^-)+\Delta g .
\end{align*} 
As in the proof of Lemma \ref{lemma:gf}, it can be shown that the modified profiles define a feasible and optimal solution. However, we also have $\sum_{t\in[T]}|\tilde r_1(t)|<\sum_{t\in[T]}|r^*_1(t)|$ and $\sum_{t\in[T]}\tilde g_1(t)=\sum_{t\in[T]}g^*_1(t),$  which contradicts our hypothesis.
%

(b) This proof is similar to the first part. Suppose to the contrary
$g^*_1(t^+)>g^*_1(t^-)\geq 0$. Define 
\begin{align*}
\Delta g': = \min\Big\{ r_1^*(t^+), -r^*_1(t^-),& g^*_1(t^+)-g^*_1(t^-), \\
&~ \min_{t\in[t^-+1,t^+-1]} (b^*_1-s^*_1(t))\Big\}.
\end{align*}
Again, we have $\Delta g'>0$ and thus we construct modified profiles $\tilde{g}_1(t), \tilde{r}_1(t)$ that differ from ${g}^*_1(t), r^*_1(t)$ only at $t^+$ and $t^-$ as follows: 
\begin{align*}
\tilde{g}_1(t^+) = g^*_1(t^+) - \Delta g', \quad \tilde{g}_1(t^-) = g^*_1(t^-) +  \Delta g', \\
\tilde{r}_1(t^+) = r^*_1(t^+) - \Delta g', \quad \tilde{r}_1(t^-) = r^*_1(t^-)+\Delta g'.
\end{align*} 
It can be checked that the modified profiles define a feasible and optimal solution with $\sum_{t\in[T]}|\tilde r_1(t)|<\sum_{t\in[T]}|r^*_1(t)|,$ and $\sum_{t\in[T]}\tilde g_1(t)=\sum_{t\in[T]}g^*_1(t)$, contradicting our hypothesis.
\end{IEEEproof}

\vspace{5pt}


Next, we choose scalars $ \beta_i$ for each $i=0,1,\ldots,m$ that helps us in defining $\tilde{r}(t)$. First, define 
\begin{align*}
t^+_0 :=\arg\max_t\left\{g^*_1(t)\ | \ t\in[\TT_0 + 1,\TT_{1}] \text{ and } r_1^*(t)>0 \right\}, 
\end{align*}
and for $i=1,\ldots,m$, 
\begin{align*}
t^-_i :=\arg\min_t\left\{g^*_1(t)\ | \ t\in[\TT_i+1,\TT_{i+1}] \text{ and } r_1^*(t)<0 \right\}. 
\end{align*}
Notice that $t_0^+$ always exists. However, $t^-_i$ may not be defined if $\TT_i + 1 = \TT_{i+1}$. Now, define the scalars $\beta_i$'s as follows.
\begin{align}
\beta_0 &:=  g^*_1(t^+_0),\label{eq:beta_0}\\ 
\beta_i &:= \begin{cases} 
g^*_1(t^-_i) & \text{ if } \TT_i+1<\TT_{i+1},\label{eq:beta_i_1}\\
1/(2\xi) &\text{ else },
 \end{cases}  \ \ i=1,\ldots,m-1, \\
 \beta_m &:=  g^*_1(t^-_m).\label{eq:beta_i_2}
\end{align}
We characterize the properties of $\beta_i$'s in the following result.
\begin{lemma}\label{lem:4}
For $i=0,1,\ldots,m$,
\begin{enumerate}[(a)]
\item
$0\leq\beta_0\leq\beta_i\leq 1/(2\xi)$,
\item $g^*_1(t)r_1^*(t) \leq \beta_ir_1^*(t)$ for all $t\in[\TT_i+1,\TT_{i+1}]$.
\end{enumerate}
\end{lemma}
\begin{IEEEproof}
(a)
From \eqref{eq:beta_0}-\eqref{eq:beta_i_2}, it is straightforward to conclude that $\beta_i$'s are nonnegative. Next, we show the following property: for any $t^-\in[\TT_0+1,\TT_{m+1}]$ 
such that $r_1^*(t^-)<0$, we have
\begin{align}\label{eq:g_suff}
g_1^*(t^-) \leq 1/(2\xi).
\end{align}
Suppose to the contrary that $g_1^*(t^-) > 1/(2\xi)$. 
Define
$$
\Delta g := \min\left\{ g_1^*(t^-) ~,~ g_1^*(t^-)-r_1^*(t^-) - 1/(2\xi) \right\}.
$$
Notice that $\Delta g>0$. 
We define a modified generation profile $\tilde{g}(t)$ that differs from $g^*(t)$ only at $t^-$. Specifically, $\tilde g_1(t^-) :=  g_1^*(t^-) - \Delta g$.
Clearly, 
$0\leq \tilde g_1(t^-) < g_1^*(t^-) \leq \overline{g}_1$.
Furthermore, from the definition of $\Delta g$, it follows that 
\begin{align}
1/(2\xi) \leq \tilde g_1(t^-) - r_1^*(t^-) &\leq  g_1^*(t^-) - r_1^*(t^-) \label{eq:tilde_p12}
\\  &= p_{12}^*(t^-)\leq f_{12}.\notag
\end{align}
Now, the function $F(x) := x-\xi x^2$ is decreasing for $x>1/(2\xi)$. Thus, we have
$$F(\tilde g_1(t^-) - r_1^*(t^-) ) \geq F(g_1^*(t^-) - r_1^*(t^-)) \geq d_2(t^-) + r_2^*(t^-),$$
i.e., $\tilde g_1$ satisfies \eqref{eq:weKnow} at $t^-$. Hence, $\tilde g_1$ is feasible. Furthermore, $\tilde g_1(t)\leq g^*_1(t)$ for all $t$. Since the cost of generation is nondecreasing, it implies that $\tilde g_1$ is in fact optimal. However, $\sum_{\tau\in[T]}\tilde g_1(\tau)<\sum_{\tau\in[T]} g^*_1(\tau)$ which contradicts our initial assumption. This completes the proof of \eqref{eq:g_suff}.

Now we turn to prove $\beta_0\leq\beta_i\leq 1/(2\xi)$. Let $t^-$ in \eqref{eq:g_suff} equal $t^-_i$. Using \eqref{eq:g_suff} in the definition of $\beta_1, \ldots, \beta_m$, it is easy to conclude that $\beta_i\leq 1/(2\xi)$, for each $i=1,\ldots,m$.

Finally, for $i = 1, \ldots, m$, if $t^-_i$ exists, then choose $t^+ = t_0^+$ and $t^-=t^-_i$. Applying Lemma \ref{lem:12}, we get $\beta_0 = g_1^*(t^+_0)\leq g_1^*(t^-_i) = \beta_i$. If $t^-_i$ does not exist, let $t^-$ be any instant after $t_0^+$, where $r_1^*(t) < 0$. Again, applying Lemma \ref{lem:12} with $t^+ = t_0^+$ and using \eqref{eq:g_suff}, we get $\beta_0 = g_1^*(t^+_0)\leq 1/(2\xi) = \beta_i$. That completes the proof of Lemma \ref{lem:4}(a).




(b) It suffices to prove that
\begin{align}\label{eq:show3}
\begin{cases} g_1^*(t)\leq \beta_i &\text{ if } r_1^*(t)>0, \\
g_1^*(t)\geq \beta_i &\text{ if } r_1^*(t)<0. \end{cases}
\end{align}
First, consider the case $t\in[\TT_0+1,\TT_{1}]$. If $r_1^*(t)>0$, then from the definition of $\beta_0$, it follows that $g_1^*(t)\leq \beta_0$. If, on the other hand, $r_1^*(t)<0$, then from Lemma \ref{lem:12}, we have $g_1^*(t)\geq g_1^*(t_0^+) = \beta_0$. 

Next, consider the case where $t\in[\TT_i+1,\TT_{i+1}]$ for some $i \in \{1,\ldots,m \}$. If $\TT_i+1=\TT_{i+1}$, then $r_1^*(\TT_{i+1})=0$ and \eqref{eq:show3} holds trivially. Otherwise, $\beta_i = g_1^*(t_i^-)$ from \eqref{eq:beta_i_1}. If $r_1^*(t)<0$, then $g_1^*(t)\geq g_1^*(t_i^-) = \beta_i$ from the definition of $t_i^-$. If $r_1^*(t) > 0$, then $g_1^*(t)\leq g_1^*(t_i^-) = \beta_i$  from Lemma \ref{lem:12}. This completes the proof of Lemma \ref{lem:4}(b).
\end{IEEEproof}
%

Now, we use the $\beta_i$'s to construct $\tilde r(t)$ as follows. For $i=0,1,\ldots,m$, and $t\in[\TT_{i}+1,\TT_{i+1}]$:
\begin{align*}
\tilde r(t) :=  (1-2\xi \beta_i) r^*_1(t) - \begin{cases}
 2\xi (\beta_m - \beta_0) b^*_1 & t=\TT_{m+1} \\ 0
 & \text{else}.
\end{cases}
\end{align*}

It is  left to prove that a storage placement solution defined by $g_1(t)=g_1^*(t)$, $r_1(t)=0$, $r_2^*+ \tilde r(t)$, $b_1 = 0$, $b_2=b_2^*+b_1^*$ is indeed feasible for $\hat P^L$. Denote $\tilde s(t) := \sum_{\tau=1}^{t}\tilde r(\tau).$ We show this in the following four steps. 


$1)~$\textit{To verify $g_1^*(t) - \xi (g^*_1(t) )^2 \geq d_2 (t) + r^*_2(t) + \tilde r(t)$:}
For any $t\in[\TT_{i}+1,\TT_{i+1}]$, it follows from \eqref{eq:weKnow} that
\begin{align*}
g^*_1(t) - \xi (g^*_1(t))^2 &\geq d_2 (t) + r^*_2(t) + \left( 1 -2\xi g_1^*(t)\right)r_1^*(t).
\end{align*}
When combined with Lemma \ref{lem:4}(b), we get
\begin{align}
g^*_1(t) - \xi (g^*_1(t))^2 
&\geq d_2 (t) + r^*_2(t) + (1-2\xi\beta_i) r^*_1(t). \label{eq:y44}
\end{align}
For $t<\TT_{m+1}$ identify $\tilde{r}(t) = (1-2\xi\beta_i) r^*_1(t)$ in \eqref{eq:y44} to establish the desired result. For $t=\TT_{m+1}$ use  $\beta_m\geq \beta_0$ from  Lemma \ref{lem:4}(a) to find that $   \tilde r(\TT_{m+1})=(1-2\xi\beta_m) r^*_1(\TT_{m+1}) - 2\xi(\beta_m-\beta_0)b_1^* \leq(1-2\xi\beta_m) r^*_1(\TT_{m+1})$. Combine this with \eqref{eq:y44} to conclude.

$2)$ \textit{To verify $g_1^*(t) \leq f_{12}$}:
This is similar to the proof of Lemma \ref{lemma:gf}(\ref{lemma.2}). Consider $t_{\max}$ to be the time instant in the cycle where $g_1^*(t)$ is maximum. Using Lemma \ref{lem:12}, it can be argued that $r_1^*(t_{max})\leq 0$. Hence, $\max g^*_1(t) \leq g^*_1(t_{max}) - r^*_1(t_{max}) = p^*_{12}(t_{max}) \leq f_{12}$.


$3)$ \textit{To verify $b_2^* + b_1^* \leq h$}: This follows from the feasibility of the problem $\hat{P}^L$.

$4)$ \textit{To verify $0 \leq \tilde s(t) \leq   b_1^*$, $t \in [\TT_0, \TT_{m+1}]$ and $\tilde s(\TT_{m+1}) = 0$}:  
First, consider the case $t\in[\TT_0+1,\TT_1]$. We have $0\leq s^*_1(t)\leq b_1^*$. Also, $\beta_0\leq 1/(2\xi)$ from Lemma \ref{lem:4}(a). Combining the above, 
\begin{align*}
0\leq \tilde s(t) = (1-2\xi\beta_0) s^*_1(t)\leq b_1^*.
\end{align*}

Next, consider the case $t\in[\TT_i+1,\TT_{i+1}]$, for $i \in \{1,\ldots,m -1\}$. Note that $s^*_1(\TT_j)=b_1^*$ for each $j=1,\ldots,i$. Thus, we have
\begin{align}
\tilde s(t) &= \sum_{j=0}^{i-1}\left[ (1-2\xi\beta_j) \sum_{\tau=\TT_j+1}^{\TT_{j+1}}r^*_1(\tau) \right] \notag\\
& \qquad +  (1-2\xi\beta_i) \sum_{\tau=\TT_i+1}^{t}r^*_1(\tau)  \notag\\
& = (1-2\xi\beta_0) b_1^* + (1-2\xi\beta_i)(s^*_1(t)-b_1^*)\label{eq:exp1}\\
&= (1-2\xi\beta_i)s^*_1(t) + 2\xi(\beta_i-\beta_0)b_1^*\label{eq:exp2}.
\end{align}
From \eqref{eq:exp1}, we find that $\tilde{s}(t)\leq b_1^*$. This follows from two observations: (i)  $0 \leq(1-2\xi\beta_i) \leq 1$ from Lemma \ref{lem:4}(a) and (ii) $s^*_1(t)-b_1^* \leq 0$ from the feasibility of $s^*_1(t)$. On the other hand, \eqref{eq:exp2} establishes that $\tilde s(t)\geq 0$. For this observe that: (i) $\beta_0\leq \beta_i\leq 1/(2\xi)$ from Lemma \ref{lem:4}(a) and (ii) $s^*_1(t) \geq 0$ from the feasibility of $s^*_1(t)$.
The argument for the case $t \in[\TT_m+1,\TT_{m+1} - 1]$ is similar. Finally, we analyze $\tilde s (\TT_{m+1})$. Using \eqref{eq:exp2}, we have
\begin{align*}
\tilde s (\TT_{m+1}) 
& = \tilde s (\TT_{m+1} - 1) + \tilde r(\TT_{m+1}) \\
&= (1-2\xi\beta_i)s^*_1(\TT_{m+1}-1) + 2\xi(\beta_m-\beta_0)b_1^*\\
& \qquad +  (1-2\xi \beta_m) r^*_1(\TT_{m+1}) - 2\xi (\beta_m - \beta_0) b^*_1 \\
& = (1-2\xi\beta_m)s^*_1(\TT_{m+1})  = 0.
\end{align*}



\subsection{Proof of $\hat \pi^{\{1\},L}_* = \pi^{\{1\},L}_*$}\label{sec:tight}
Among the optimal points of $\hat\Pi^{\{1\},L}$, consider the one that minimizes
$\sum_{t\in [T]}(g_1(t) + |r_2(t)|)$; let $g_1^*(t), r_2^*(t),t\in[T]$ be the corresponding optimal generation profile and storage charging/ discharging profile at node 2. If $g_1^*(t) - \xi(g_1^*(t))^2 = d_2(t) + r_2^*(t)$ for all $t\in[T]$, then $g_1^*(t),r_2^*(t)$ are feasible for $\Pi^{\{1\},L}$ and hence $\hat \pi^{\{1\},L}_* = \pi^{\{1\},L}_*$. Now, suppose to the contrary there exists $t_0\in[T]$ such that 
\begin{align}\label{eq:contra37}
g_1^*(t_0) - \xi(g_1^*(t_0))^2 > d_2(t_0) + r_2^*(t_0).
\end{align}
We consider the following two cases.
\subsubsection*{Case 1: $g_1^*(t_0) = 0$} Then $d_2(t_0) + r^*_2(t_0) < 0$ and hence $r_2^*(t_0)<0$. Let $t_1 := \max\{ t\in[1, t_0-1] \ | \ r_2^*(t_1)>0 \}$. Notice that $t_1$ always exists. Define
$$
\Delta r := \min\{ r_2^*(t_1), -(d_2(t_0) + r_2^*(t_0)), \min_{t\in[t_1+1,t_0-1]} s_2^*(t) \}.
$$
It can be checked that $\Delta r>0$. Construct modified storage profile $\tilde r_2(t),t\in [T]$ that differs from $r_2^*(t),t\in [T]$ only at $t_1$ and $t_0$ as follows:
\begin{align*}
\tilde r_2(t_1) := r_2^*(t_1) - \Delta r \quad \text{ and } \quad
\tilde r_2(t_0) :=  r_2^*(t_0) + \Delta r.
\end{align*}
Then, it can be shown that $0\leq \sum_{\tau=1}^{t}{\tilde r_2(\tau)}\leq b_2, t\in[T]$ and $\sum_{\tau\in[T]}{\tilde r_2(\tau)}=0$.
Also, $g^*(t_1)-\xi g^*(t_1)\geq d_2(t_1) +r^*_2(t_1)\geq d_2(t_1) + \tilde r_2(t_1)$ and  $d_2(t_0) + \tilde r_2(t_0) \leq 0 =  g^*_1(t_0)$. Thus, the solution defined by the modified profile is feasible. It is also optimal since it has the same generation profile $g_1^*(t),t\in[T]$. However, $\sum_{t\in[T]}|\tilde r_2(t)| < \sum_{t\in[T]}| r^*_2(t)|$ which leads to a contradiction.

\subsubsection*{Case 2: $g_1^*(t_0) > 0$} 
Here, we construct a modified generation profile  $\tilde g_1(t),t\in [T]$ as follows. If $d_2(t_0) + r_2^*(t_0) < 0$, let $\tilde g_1(t_0) = 0$. Otherwise, set $\tilde g_1(t_0)$ such that $\tilde g_1(t_0) - \xi (\tilde g_1(t_0))^2 = d_2(t_0) + r_2^*(t_0)$ and $\tilde g_1(t_0) \leq 1 / (2\xi)$ . It is easy to check that such a $\tilde g_1(t_0)$ exists since the quadratic function $F(x) := x - \xi x^2$ increases in $[0, 1/(2\xi))$ and decreases in $[1/(2\xi), 1/\xi]$. In both cases, $\tilde g_1(t_0) < g_1^*(t_0)$ and hence $\tilde g_1$ defines an optimal generation profile with $\sum_{\tau\in [T]}\tilde g_1(\tau)<\sum_{\tau\in [T]}\tilde g^*_1(\tau)$, which is a contradiction.


\begin{figure*}[t]
	\subfloat[]{
        		\centering
		\includegraphics[width=0.31\textwidth]{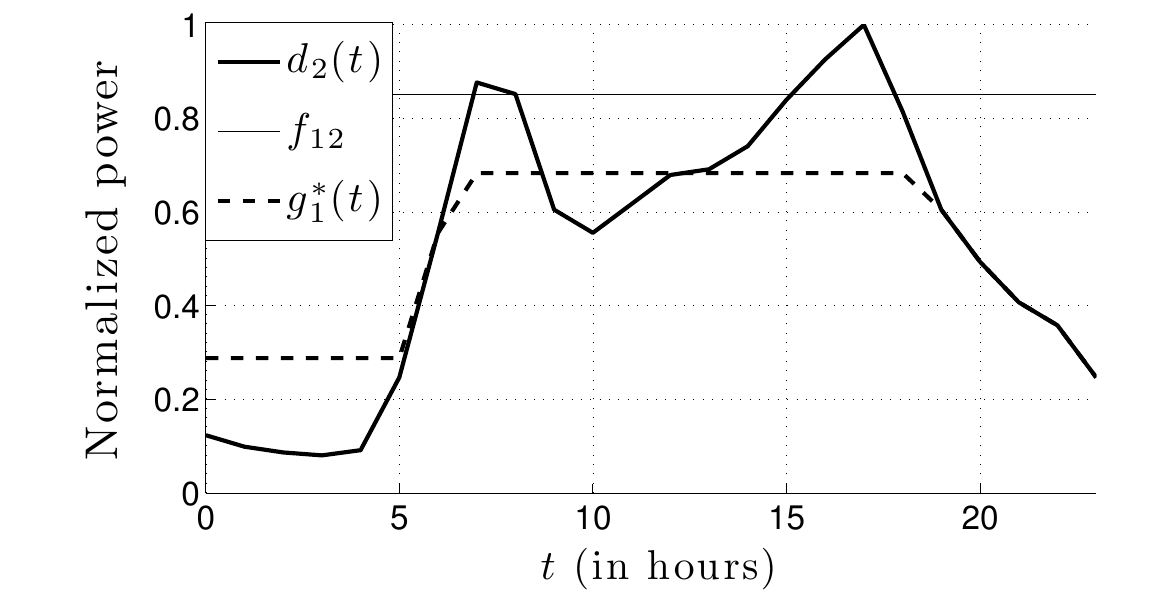}
		\label{fig:sample}
	}
	\subfloat[]{
                \centering
		\includegraphics[width=0.31\textwidth]{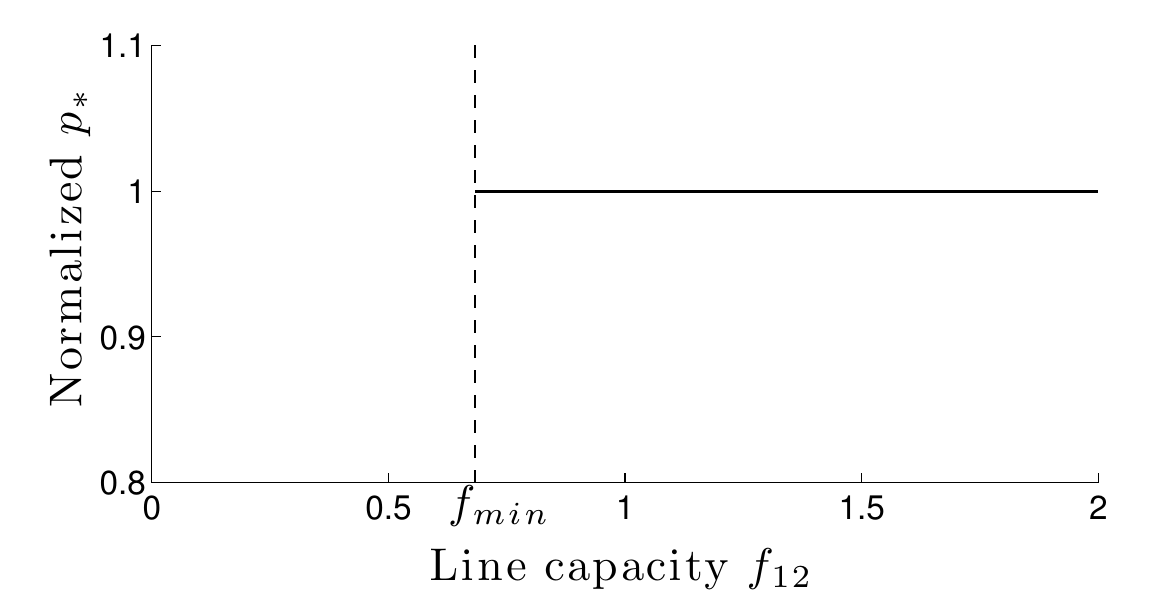}
                \label{fig:sample.f}
         }
        \subfloat[]{
                \centering
		\includegraphics[width=0.31\textwidth]{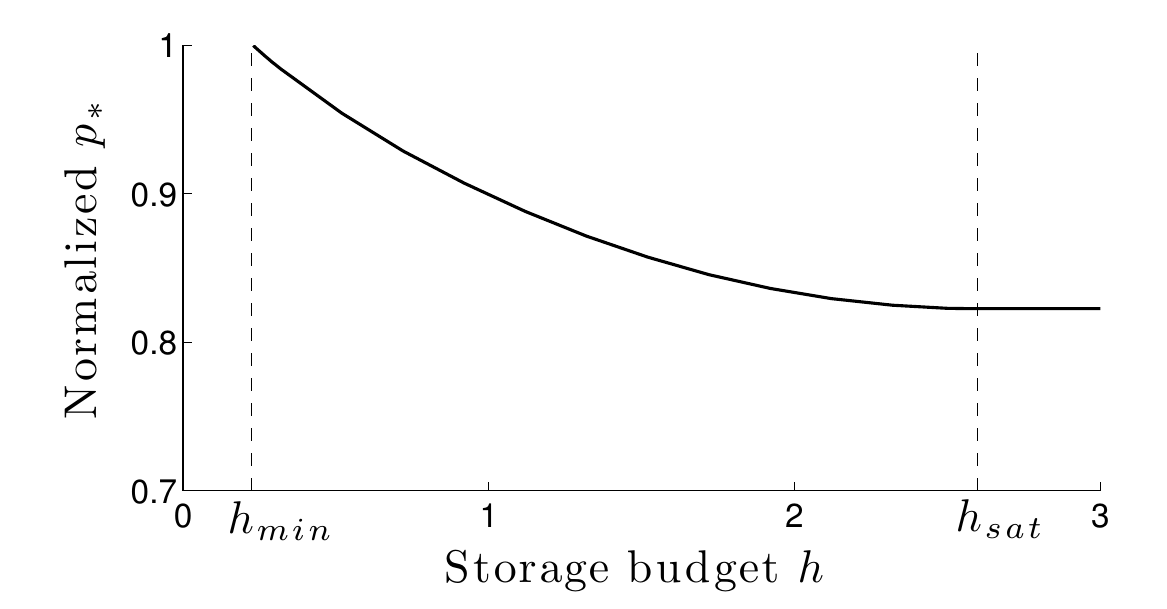}
                \label{fig:sample.h}
	}
\caption{{Plots to illustrate Propositions \ref{prop:SGSLf} and \ref{prop:SGSLhmin} in the Appendix. All quantities are in per units.  (a) Typical hourly load profile $d_2(t), t\in [T]$ and optimal generation portfolio for line flow capacity $f_{12}=0.85$, generation capacity $\overline{g}_1=1$ and  storage budget $h=1$. Notice that $\max_{t\in [T]} g_1^*(t) \leq f_{12}$ as stated in Lemma \ref{lemma:gf}.
 (b) $p_*(\overline{g}_1=1,f_{12}=0.85,h)$, for $h$ in $[0, 3]$. Notice that $f_{12} \leq \max_{t\in [T]} d_2(t)$, i.e., the problem is infeasible in the absence of storage. $h_{\min} = 0.226$ and $h_{sat} = 2.598$ were calculated using Proposition \ref{prop:SGSLhmin}. (c) $p_*(\overline{g}_1=1,f_{12},h=1)$  for $f_{12}$ in $[0, 2]$. As in Proposition \ref{prop:SGSLf}, the problem is infeasible for $f_{12}<f_{\min} = 0.683$ and the optimal cost remains constant for $f_{12}\geq f_{\min}$.}}
	\label{fig:sgslSims}
\end{figure*}

\section{Partial results on specific network topologies}
\label{app:partial}
Here, we present some partial results on the storage placement problems for two specific classes of networks: (1) a 2-node network with a single generator bus and a single load bus in Figure \ref{fig:SGSL}, and (2) a star network with a central generator bus serving multiple loads via transmission-constrained lines. We explore some salient features of the optimal placement. For ease of presentation, we assume that the costs $c_k(\cdot)$ are strictly convex at each node, neglect storage efficiency losses $\ac = \ad = 1$ and ignore ramping constraints with $\epsilon_\gamma=\epsilon_\delta=1$. The proofs of the results are included in Appendix \ref{app:proofs}.

\subsection{Single Generator Single Load Network}
\label{sec:SGSL}

For this network, placing all the available storage resources at the load bus is always optimal. This is an immediate consequence of Theorem \ref{thm:main}. Now, for any fixed demand profile $d_2(t),t \in [T]$ at the load bus, we analyze the behavior of the optimal cost of production as a function of the generation capacity $\overline{g}_1$, the line flow capacity $f_{12}$ and the available storage budget $h$. Through this analysis, we hope to gain insights on marginal savings in terms of generation cost with additional investments in storage capacity, line capacity and storage budget.
Let the parameterized storage placement problem be $P(\gc, f_{12}, h)$ and its optimal cost be $p_* (\gc, f_{12}, h)$. Similarly define, $\Pi^{\{1\}}(\gc, f_{12}, h)$ and $\pi^{\{1\}}_*(\gc, f_{12}, h)$. Then the optimal costs satisfy the following:

\begin{proposition}
\label{prop:SGSLf}
For any $h\geq 0$, problem $P(\gc,f_{12}, h)$ is feasible iff $\mingcf \geq f_{\min}$, where
\begin{align}\label{eq:fmin}
f_{\min} &= 
\max \left\{     \max_{1\leq t \leq T} \left( \frac{ \sum_{\tau=1}^{t} d_2(\tau)} {t} \right), \right. \notag \\
& \left. \qquad \qquad
\max_{1\leq t_1 < t_2 \leq T} \left( \frac{\sum_{\tau=t_1+1}^{t_2} {d_2(\tau)}-h}{t_2-t_1} \right)
\right\}.
\end{align}
Moreover, if $\mingcf \geq f_{\min}$, then $p_*(\gc,f_{12}, h) = p_*(f_{\min},f_{\min}, h)$.
\end{proposition}

We interpret this result as follows. If either the line flow limit $f_{12} < f_{\min}$ or the generation capacity $\gc < f_{\min}$, the load cannot be satisfied. Notice that $f_{min}$ is non-increasing in $h$; this suggests that storage can play a role in avoiding transmission upgrades when dealing with lowering generation costs \cite{denholm09}. However, when $f_{12} \geq f_{\min}$ and $\gc \geq f_{\min}$, the optimal cost of operation does not depend on the specific values of $f_{12}$ and $\gc$. Thus, investing in line limits and generation capacities over $f_{\min}$  do not reduce the cost of operation. The behavior of $p^*(\gc,f_{12}, h)$ as a function of $f_{12}$ is illustrated in Figure \ref{fig:sample.f}.

Next, we characterize the behavior of $P(\gc,f_{12}, h)$ and its optimal cost $p^*(\gc,f_{12}, h)$ as a function of $h$. For a given $f_{12}$ and $\overline{g}_1$, the minimum required storage budget to serve the load depends on the demand profile $d_2(t), t\in[T]$; it may even be the case that the problem remains infeasible no matter how large the storage budget is. Also, when we allow larger storage budgets, the generation cost does not reduce beyond a point, i.e., there exists $h_{sat}$ such that $p_*(\gc,f_{12}, h) = p_*(\gc,f_{12}, h_{sat})$ for all $h \geq h_{sat}$. Between $h_{min}$ and $h_{sat}$, the optimal cost $p^*(\gc,f_{12}, h)$ is decreasing and convex in $h$. Thus, we have diminishing marginal returns on the investment on storage, i.e., the benefit of the first unit installed is more than that from the second unit. The behavior of $p^*(\gc,f_{12}, h)$ is illustrated in Figure \ref{fig:sample.h}. In what follows, we formalize this discussion.

Construct the sequence $\{\tau_m\}_{m=0}^M$ as follows. Let $\tau_0 = 0$. Define $\tau_m$ iteratively:
\begin{align}
\label{eq:tn}
\tau_m & = \argmax_{\tau_{m-1}+1\le t\le T} {\left( {\frac{\sum_{\tau=\tau_{m-1}+1}^{t} {d_2(\tau)}}{t-\tau_{m-1}} } \right) },
\end{align}
for $1\leq m \leq M$, where $M$ is the smallest integer for which $\tau_M=T$. Note that the sequence depends only on the demand profile $d_2(t), t \in [T]$. For any $x \in \mathbb{R}$, let $[x]^+ := \max(x, 0)$. 

\begin{proposition}\label{prop:SGSLhmin}
Problem $P(\gc,f_{12}, h)$ satisfies:
\begin{enumerate}[(a)]
\item  If $\mingcf < {\max}_{t \in [T]} \left({\frac{\sum_{\tau=1}^{t} d_2(\tau) } {t} } \right) $, then $P(\gc,f_{12}, h)$ is infeasible for all $h\geq 0$.
\item Suppose, $\mingcf \geq {\max}_{t \in [T]} \left({\frac{\sum_{\tau=1}^{t} d_2(\tau) } {t} } \right) $. Then, $P(\gc,f_{12}, h)$ is feasible iff $h \geq h_{\min}$ and $p_*(\gc,f_{12}, h)$ is convex and non-increasing in $h$, where
\begin{align}\label{eq:hmin}
h_{min}= 
\max_{0\leq t_1\leq t_2\leq T} {\left[\sum_{\tau=t_1+1}^{t_2} {(d_2(\tau)-\mingcf)}\right]^+}.
\end{align}
Furthermore, $p_*(\gc,f_{12}, h)$ is constant for all $h\geq h_{sat}$, where
\begin{align} 
\label{eq:h_o}
h_{sat} = \max_{1\leq m \leq M} u_m,
\end{align}
where
\begin{align}
u_m &:= 
\max_{\tau_{m-1}+1\le t \le \tau_m} \left[      \left(\sum_{\tau=\tau_{m-1}+1}^{\tau_m} {d_2(\tau)}\right)  {\frac{t-\tau_{m-1}}{\tau_m-\tau_{m-1}} } \right. \notag\\
& \left. \qquad \qquad \qquad \qquad \qquad - \left(\sum_{\tau=\tau_{m-1}+1}^{t}{d_2(\tau)} \right) 
\right].
\end{align}
\end{enumerate}
\end{proposition}

\subsection{Star Network}
\label{starNetwork}

\begin{figure}[h!]
	\centering
	\includegraphics[width=0.11\textwidth]{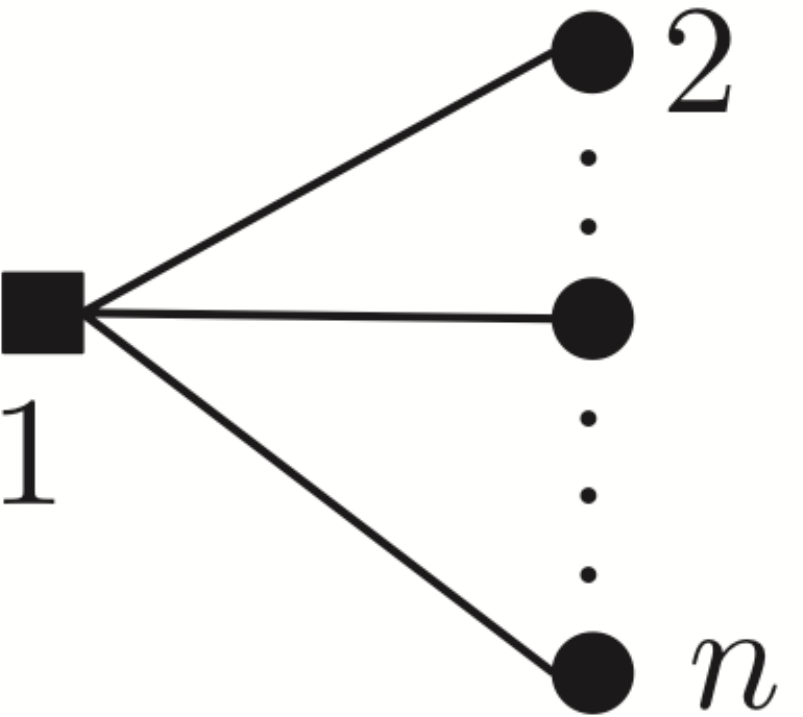}
	\caption{ An $n$-node network with a central generator bus serving $n-1$ load buses via transmission-constrained links.}
\label{fig:starNet}

\end{figure}

Consider the star network on $n\geq 2$ nodes as shown in Figure \ref{fig:starNet}. For this network $\mathcal{N}_G = \{ 1 \}$ and $\mathcal{N}_D = \{ 2, 3, \dots, n \}$. We showed in Section \ref{sec:multGen} that for such a generator bus with multiple connections, placing zero storage at bus 1 is not necessarily always optimal for an arbitrary storage budget $h$. However, Proposition \ref{prop:hmin_starN} below shows that $b_1^*=0$ provided that the available storage budget $h$ is either (a) exactly equal to the minimum value required for which $P$ and $\Pi^{\{1\}}$ are feasible or (b) is large enough. 

To formalize notation, consider the problems $P(h)$, $\Pi^{\{1\}}(h)$ and their optimal costs $p_* (h)$ and $\pi^{\{1\}}_* (h)$ as a function of $h$. For simplicity, we consider the case $\overline{g}=\infty$. Then we have the following result.

\begin{proposition}
\label{prop:hmin_starN}
Suppose $f_{1k} \geq \max_{t\in[T]} \left(\dfrac{\sum_{\tau=1}^{t}{d_{k}(\tau)}}{t} \right)$ for all $k \in \mathcal{N}_D$. There exists $h_{min}\geq 0$, such that $P(h)$ and $\Pi^{\{1\}}(h)$ are feasible if and only if $h\geq h_{min}$, where 
\begin{equation}\label{eq:h_minRadial}
h_{min} = \sum_{k \in \mathcal{N}_D}{ 
\max_{0\leq t_1<t_2\leq T}{
\left[ \sum_{\tau=t_1+1}^{t_2}{( d_{k}(\tau)-f_{1k} )} \right]^+
} }.
\end{equation}
Also, $p_*$ and $\pi^{\{1\}}_*$ satisfy
\begin{enumerate}[(a)]
\item $p_*(h_{min})=\pi^{\{1\}}_*(h_{min})$,
\item there exists $h_{sat}\geq h_{min}$ such that $p_*(h) = \pi^{\{1\}}_*(h) $ for all $h\geq h_{sat}$.
\end{enumerate}
\end{proposition}
The behavior of $p_* (h)$ and $\pi^{\{1\}}_* (h)$ is depicted in Figure \ref{fig:counter2} for the case where $n=3$. We end this section with a remark on $h_{sat}$. The proof in Appendix \ref{app:proofs}
 essentially shows that $p_*(\infty)=\pi^{\{1\}}_*(\infty)$, and $h_{sat} := \max_{t \in [T]} s_2^*(t)$ where $s_2^*(t)$ is the finite optimal storage level at bus 2 for the problem $\Pi(\infty)$.

\begin{figure}[t]
     \centering
\includegraphics[scale=0.4]{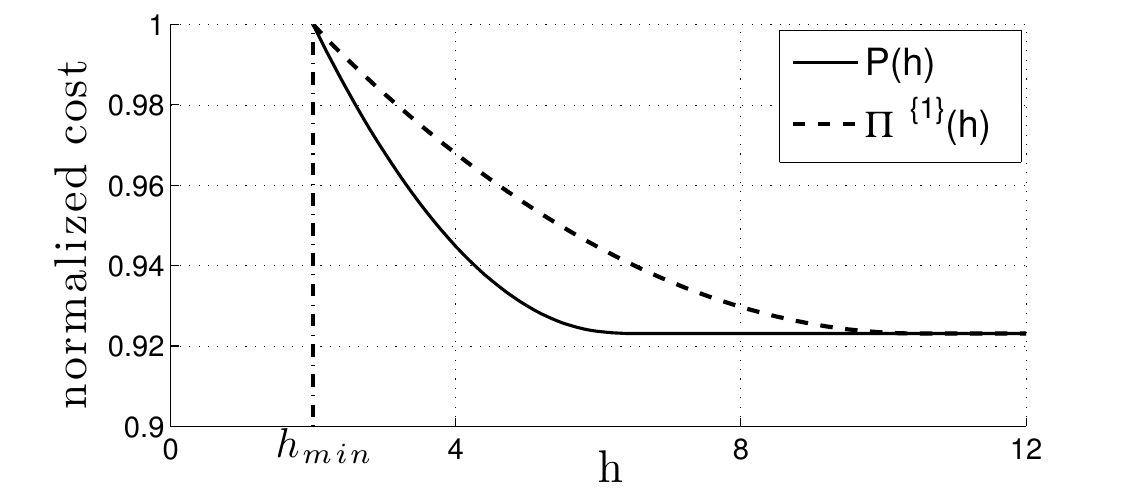}	
		\caption{$P(h)$ and $\Pi^{\{1\}}(h)$ for the simple 3-node star network ($n=3$). }
		\label{fig:counter2}
\end{figure}

Now, we restrict attention to a special case of the problem over the star network.
\begin{proposition}
\label{prop:proportionalStar}
Suppose there exists positive constants $\phi_k, k\in \mathcal{N}_D$, such that $d_k(t) = \phi_k d_2(t)$ for all $t\in[T]$ and  
\begin{align}\label{eq:pro5}
f_{1k} \geq 
\max_{0\leq t_1 < t_2 \leq T} \left( \frac{\sum_{\tau=t_1+1}^{t_2} {d_k(\tau)}}{t_2-t_1} \right).
\end{align}
If $P$ is feasible, then $\Pi^{\{1\}}$ is feasible and $p_* =\pi^{\{1\}}$.
\end{proposition}
We showed through an example in Section \ref{sec:multGen} that generator nodes with multiple connections may need nonzero storage capacities under optimal allocation. However, the demand profiles considered in the specific example showed opposite trends. In practice, aggregate demand profiles at various nodes show similar diurnal variations \cite{pjm}. To study this, we consider the case where loads at various demand nodes are proportional to each other and the line capacities are proportional to the average demands at each node. In Proposition \ref{prop:proportionalStar}, we recover a result similar to Theorem \ref{thm:main}. This suggests a potential direction for future investigation: to characterize the properties of storage placement when demands show similar variations and the transmission lines have suitable capacities to meet the demand.

\section{Proofs of partial results}
\label{app:proofs}

Here, we present the proofs of Propositions \ref{prop:SGSLf}--\ref{prop:proportionalStar}.
For the single generator single node and the star networks, we drop the voltage angles $\theta_k (t), k\in \mathcal{N}, t \in [T]$. For any value of the power flow $p_{1k}(t)$ from bus $1$ to bus $k$, voltage angles $\theta_k(t)$ can always be chosen to satisfy the power flow constraint in \eqref{eq:line flow eqn}.
Furthermore, since $\alpha=1$, define $r_k(t) := \gamma_k(t) - \delta_k(t)$ be the power that flows into the storage device at node $k\in\mathcal{N}$ at time $t\in[T]$.  Notice that $r_k(t)$ can be positive or negative depending on whether power flows in or out of the
storage device. Also, the storage level of the storage device at node $k\in\mathcal{N}$, at time $t$ can be written as $s_k(t) = \sum_{\tau = 1}^{t}{r_k(\tau)}$.



\subsection{Proof of Proposition \ref{prop:SGSLf}}\label{sec:D_1}
We drop subscripts from the variables $d_2(t), g_1(t), t \in [T]$, $f_{12}, \overline{g}_1, c_1(\cdot)$ and the superscript from $\Pi^{\{1\}}(\cdot), \pi^{\{1\}}_*(\cdot)$ for ease of notation throughout Sections \ref{sec:D_1} and \ref{sec:D_2}.


From Theorem \ref{thm:main}, it suffices to show the claim for $\Pi (\overline{g}, f, h)$ and $\pi_* (\overline{g}, f, h)$. First, we show that if $\Pi (\overline{g}, f, h)$ is feasible, then $\mingf\geq f_{min}$. Fix any $h\geq 0$ and let $g(t), t \in [T]$ be a feasible generation profile. Since $\sum_{\tau=1}^{t} r_2(\tau) = s_2(t) \geq 0$, we have for any $t\in[T]$
\begin{align}
\label{eq:fmin.1}
\max_{ t' \in [T]} g(t') 
\geq  \frac{\sum_{\tau=1}^{t} g(\tau)}{t} 
= \frac{\sum_{\tau=1}^{t} (d(\tau) +  r_2(\tau)}{t})
\geq \frac{\sum_{\tau=1}^{t}{d(\tau)}}{t}.
\end{align} 
Furthermore, for any $1\leq t_1<t_2\leq T$, the power extracted from the storage device between time instants $t_1$ and $t_2$ cannot exceed the total storage budget $h$ and hence we have \begin{align}
\max_{t'\in [T]}g(t')
\geq \frac{\sum_{\tau=t_1 + 1}^{t_2} g(\tau)}{t_2-t_1} 
&= \frac{\sum_{\tau=t_1+1}^{t_2} d(\tau) + \sum_{\tau=t_1+1}^{t_2} r_2(\tau)}{t_2-t_1}\notag\\
&\geq \frac{\sum_{\tau=t_1+1}^{t_2}{d(\tau)-h}}{t_2-t_1}.\label{eq:fmin.2}
\end{align}
Since $g(t), t \in [T]$ is feasible, $g(t) \leq \mingf$ for all $t \in [T]$. Hence, combining \eqref{eq:fmin.1} and  \eqref{eq:fmin.2}, we get
$$
\mingf \geq \max_{t'\in [T]}g(t')\geq f_{min}.
$$

Next, we show that  $\mingf\geq f_{\min}$ is sufficient for  $\Pi(\overline{g},f,h)$ to be feasible.
Consider the optimal generation profile $g^*(t), t \in [T]$ for the relaxed problem $\Pi(+\infty, +\infty, h)$.  Suppose it satisfies
\begin{equation}\label{eq:toProve}
\max_{t \in [T]} g^*(t) \leq f_{\min}. 
\end{equation}
Then $g^*(t), t\in [T]$ is also feasible and optimal for problem $\Pi(\overline{g}, f, h)$ for $\mingf \geq f_{min}$. Also, $\pi_*(\overline{g},f,h) = \pi_*(f_{\min},f_{\min},h)$ for $\mingf \geq f_{min}$. It remains to show that \eqref{eq:toProve} indeed holds. Consider the following notation.
\begin{align*}
t_{\max} &:= \max \{ t  \in [T]\ | \ g^*(t) = \max_{\tau \in [T]} g^* (\tau) \}, \\
t_{less} &:= \max \left\{ 0 \leq t < t_{max} \ | \ g^*(t) < g^*(t_{max}) \right\}.
\end{align*}
In the above definition $g^*(0):=0$ for convenience. 
If $g^*(t_{max}) = 0$, then \eqref{eq:toProve} clearly holds. Henceforth, assume $g^*(t_{max})>0$.  Then, $g^*(t),t\in[T]$ satisfies:
\begin{align*}
\max_{t\in[T]}{g^*(t)}
&= \frac{\sum_{\tau=t_{less}+1}^{t_{max}}{g^*(\tau)}}{t_{max} - t_{less}} = \frac{\sum_{\tau=t_{less}+1}^{t_{max}}{[d(\tau)+r_2^*(\tau)]}}{t_{max} - t_{less}}
\end{align*}
\begin{align}\label{eq:check1}
= \frac{1}{t_{max} - t_{less}} \left[\left(\sum_{\tau=t_{less}+1}^{t_{max}}{ d(\tau) } \right) +
{s_2^*(t_{max})}-
{s_2^*(t_{less})} \right]
\end{align}
Now, suppose the following holds:
\begin{align}\label{eq:s2(tmax)}
s_2^*(t_{max}) = 0 \quad \text{and} \quad s_2^* (t_{less}) &= \begin{cases} 0, & \text{ if } t_{less} = 0, \\ h, & \text{ otherwise}.\end{cases}
\end{align}
If \eqref{eq:s2(tmax)} holds, it follows from \eqref{eq:check1}:
\begin{align*}
\max_{t\in[T]}{g^*(t)}  &=\begin{cases} 
\dfrac{\sum_{\tau=1}^{t_{max}}{d(\tau)}}{t_{max}},  &\text{ if } t_{less}=0,\\
\dfrac{\sum_{\tau=t_{less}+1}^{t_{max}}{d(\tau)-h}}{t_{max}-t_{less}},  & \text{ otherwise},
\end{cases} \\
&\leq f_{\min}.
\end{align*}
and hence \eqref{eq:toProve} is satisfied. Next, we show that \eqref{eq:s2(tmax)} indeed holds to complete the proof.  First we prove that $s^*_2(t_{max})=0$, i.e., the storage device at node $2$ fully discharges at time $t_{max}$. Suppose $s^*_2(t_{max})>0$. As in Lemma \ref{lemma:gf}, we construct a modified generation profile and storage control policy that is feasible and has an objective function value no greater than $\pi_*(+\infty, +\infty, h)$. But, the optimal generation profile $g^*(t), t\in[T]$ is unique since the cost function $c(\cdot)$ is assumed to be strictly convex. Hence we derive a contradiction. By hypothesis, $s^*_2(t_{max})>0$ and hence storage device at bus $2$ discharges for some $t>t_{max}$. Let $t_1$ be the first such time instant. Define
$$ \Delta_1 := \min\left\{ s_2^*(t_{max}) \ , \ g^*(t_{max}) \ , \ g^*(t_{max}) - g^*(t_1)  \right\}.$$
Notice that $\Delta_1>0$. Consider the modified generation profile $\tilde{g}(t)$ and control policy $\tilde{r}_2(t)$, that differ from $g^*(t)$ and $r_2^*(t)$ only at $t_{max}$ and $t_1$ as follows:
\begin{align*}
\tilde{g}(t_{max}) = g^*(t_{max}) - \Delta_1, \quad \tilde{g}(t_1) = g^*(t_1) + \Delta_1, \\
\tilde{r_2}(t_{max}) = r_2^*(t_{max}) - \Delta_1, \quad \tilde{r_2}(t_1) = r_2^*(t_1) + \Delta_1.
\end{align*}
Using Lemma \ref{lemma:meanVal}, we have
$$
c(\tilde{g}(t_{max})) + c(\tilde{g}(t_1)) \leq c({g}^*(t_{max})) + c({g}^*(t_{1})).
$$
It can be checked that the modified profiles are feasible for $\Pi(+\infty, +\infty, h)$. The details are omitted for brevity. This is a contradiction and hence $s^*_2(t_{max})=0$.

Next, we characterize $s_2^*(t_{less})$. If $t_{less}=0$, then $s_2^*(t_{less})=s_2^0 = 0$. If $t_{less}>0$, we prove that $s_2^*(t_{less}) = h$, i.e., the storage device at node $2$ is fully charged at time $t_{less}$. Suppose $s_2^*(t_{less}) < h$. As above, we construct a modified generation profile $\tilde{g}(t)$ and storage control policy $\tilde{r}_2(t)$ that achieves an objective value no greater than $\pi_*(+\infty, +\infty, h)$ to derive a contradiction. In particular, define
$$ \Delta_2 := \min\left\{ h - s_2^*(t_{less}) , g^*(t_{less}+1) , g^*(t_{less}+1) - g^*(t_{less})  \right\}.$$
Observe, $\Delta_2>0.$
Consider $\tilde{g}(t)$ and $\tilde{r}_2(t)$, that differ from $g^*(t)$ and $r_2^*(t)$ only at $t_{less}$ and $t_{less}+1$ as follows:
\begin{align*}
\tilde{g}(t_{less}) = g^*(t_{less}) + \Delta_2, \quad \tilde{g}(t_{less}+1) = g^*(t_{less}+1) - \Delta_2, \\
\tilde{r_2}(t_{less}) = r_2^*(t_{less}) + \Delta_2, \quad \tilde{r_2}(t_{less}+1) = r_2^*(t_{less}+1) - \Delta_2.
\end{align*}
As above, this defines a feasible point for $\Pi(+\infty, +\infty, h)$ and achieves an objective value strictly less  than $\pi_*(+\infty, +\infty, h)$. This is a contradiction and hence $s_2^*(t_{less}) = h$ for $t_{less} > 0$.


\subsection{Proof of proposition \ref{prop:SGSLhmin}}\label{sec:D_2}

From Theorem \ref{thm:main}, it suffices to prove the claim for $\Pi(\overline{g}, f, h)$ and $\pi_*(\overline{g}, f, h)$.

${(a)}$
To the contrary of the statement of the Proposition suppose that  $\mingf < {\max}_{t \in [T]} \left({\frac{\sum_{\tau=1}^{t} d(\tau) } {t} } \right) $ and $\Pi(\overline{g}, f, h)$ is feasible for some $h\geq 0$.  Then, it  follows directly from Proposition \ref{prop:SGSLf} that $\mingf \geq f_{min} \geq {\max}_{t \in [T]} \left({\frac{\sum_{\tau=1}^{t} d(\tau) } {t} } \right)$, contradicting our hypothesis.

${(b)}$
First we show that if $\Pi(\overline{g}, f, h)$ is feasible then $h\geq h_{min}$. Suppose $\Pi(\overline{g}, f, h)$ is feasible. Then, for all $0 \leq t_1 < t_2 \leq T$ Proposition \ref{prop:SGSLf} implies that $\mingf \geq f_{min}\geq   \left( \sum_{\tau=t_1+1}^{t_2}{d(\tau)}-h \right)/(t_2-t_1).$ Rearranging this we get $
h \geq \sum_{\tau=t_1+1}^{t_2} \left(  d(\tau)  - \mingf \right).
$ Also, $h \geq 0$ and hence:
\begin{align*}
 h & \geq   \max_{0\leq t_1< t_2\leq T} {\left[\sum_{\tau=t_1 + 1}^{t_2} {(d(\tau) - \mingf)}\right]^+ } =  h_{min} \notag.
\end{align*}

Now we show that  $h \geq h_{min}$ is sufficient for $\Pi(\overline{g}, f, h)$ to be feasible. The relation $h \geq h_{min}$ can be equivalently written as follows:
\begin{align}
\label{eq:hhmin.2}
\mingf \geq \dfrac{\sum_{\tau=t_1+1}^{t_2}{d(\tau)}-h}{t_2-t_1},\quad \text{for all } 0\leq t_1 < t_2 \leq T.
\end{align}
Also, by hypothesis, we have 
\begin{align}
\label{eq:hhmin.3}
\mingf \geq \displaystyle\max_{t \in [T]} \left({\frac{\sum_{\tau=1}^{t} d(\tau) } {t} } \right).
\end{align}
Combining \eqref{eq:hhmin.2} and \eqref{eq:hhmin.3}, we get $\mingf \geq f_{min}$. Then, Proposition \ref{prop:SGSLf} implies that $\Pi(\overline{g}, f, h)$ is feasible. Convexity and non-decreasing nature of $p_*(\overline{g},f,h)$ as a function of $h$ follows from linear parametric optimization theory \cite{boyd_book}.

Finally, we prove that $p_*(\overline{g},f,h)$ is constant for all $h\geq h_{sat}$, where $h_{sat}$ is as defined in \eqref{eq:h_o}.
The proof idea here is as follows. We construct the optimal generation profile $g^*(t), t \in [T]$ for the problem $\Pi(+\infty, +\infty, +\infty)$ and show that it is feasible and hence optimal for the problem $\Pi( \overline{g}, f, +\infty)$ provided $\mingf \geq {\max}_{ t\in[T]} \left({\frac{\sum_{\tau=1}^{t} d(\tau) } {t} } \right)$ holds. Problem $\Pi(+\infty, +\infty, +\infty)$ can be re-written as follows.
\begin{subequations}
\begin{align}
\underset{g(t), t\in [T]}{\text{minimize}}  & \ \ \ \ \sum_{t =1}^{T}{ c_1 \left(g(t)\right) } \notag\\
\text{subject to} 
& \ \ \ g(t) \geq 0, \quad \sum_{\tau = 1}^t ( g(\tau) - d(\tau) ) \geq 0, \qquad t \in [T], \label{eq:pi.2}\\
& \ \ \ \sum_{\tau = 1}^T g(\tau) = \sum_{\tau = 1}^T d(\tau). \label{eq:pi.3}
\end{align}
\end{subequations}
Let the Lagrange multipliers in equations \eqref{eq:pi.2}--\eqref{eq:pi.3} be $\lambda(t)$, $\ell(t)$, $t \in [T]$ and $\nu$, respectively. 
It can be checked that the following primal-dual pair satisfies the Karush-Kuhn-Tucker conditions and hence is optimal for the convex program $\Pi(+\infty, +\infty, +\infty)$ and its Lagrangian dual \cite{boyd_book}. We omit the details for brevity.
\begin{subequations}
\begin{align*}
g^*(t) &= \dfrac{\sum_{\tau=\tau_{m-1}+1}^{\tau_m}{d(\tau)} }{ \tau_m-\tau_{m-1} } , \ \  t=[\tau_{m-1}+1,\tau_m],\   m=1,\dots,M,  \\
\ell^*(t) &= \begin{cases}
        c'(g^*(\tau_m))-c'(g^*(\tau_{m}+1)), & \text{ if } t=\tau_{1},\tau_{2},\ldots,\tau_{M-1}\\
        0, & \text{otherwise}
        \end{cases}\\
\lambda^*(t) &= 0 , \quad t\in[T], \quad \text{ and } \quad \nu^* = -c'(g^*(T)) .
\end{align*}
\end{subequations}
The above profile $g^*(t), t \in [T]$ of $\Pi(+\infty, +\infty, +\infty)$ satisfies:
\begin{align*}
\max_{t \in [T]}{g^*(t)} = \max_{t \in [T]}{\frac{\sum_{\tau=1}^{t}{d(\tau)}}{t}} \leq \mingf,
\end{align*}
and hence is feasible and optimal for $\Pi(\overline{g}, f, +\infty)$. 
Note that $\sum_{\tau=\tau_{m-1}+1}^{\tau_{m}}{\left(g^*(\tau)-d(\tau)\right)}=0$ for all $1\leq m\leq M$. Thus, for $\tau_{m-1}< t \leq \tau_m$, we have 

\begin{align*}
s_2^*(t) &= \sum_{\tau = \tau_{m-1}}^{t}{\left(g^*(\tau)-d(\tau)\right)} \\
&=  \frac{\sum_{\tau=\tau_{m-1}+1}^{\tau_m} {d(\tau)}}{\tau_m-\tau_{m-1}} (t-\tau_{m-1}) - \sum_{\tau=\tau_{m-1}+1}^{t}{d(\tau)}.
\end{align*}
Maximizing the above relation over all $t\in[T]$ we get $\max_{t\in[T]}{s_2^*(t)} = h_{sat}$.  Therefore, $g^*(t),t\in [T]$ is feasible and optimal for $\Pi(\overline{g}, f, h)$ provided that $h\geq h_{sat}$. 

\subsection{Proof of Proposition \ref{prop:hmin_starN}}
%

%
First we show that $h\geq h_{min}$ is necessary for $P(h)$ to be feasible. Consider any feasible solution of $P(h)$.
For any $k\in \mathcal{N}_D$ and $0\leq t_1<t_2\leq T$, we have 
$
\sum_{\tau=t_1+1}^{t_2}{r_{k}(\tau)}\geq -b_k,
$
since the power extracted from the storage device at node $k$ cannot exceed the corresponding storage capacity $b_k$. Also, for any $k\in\mathcal{N}_D$ the power flow on the line joining buses 1 and $k$ satisfies $p_{1k}(t) = d_k(t) + r_k(t)\leq f_{1k}$ for all $t\in [T]$. Combining the above relations and rearranging, we get
$
b_k \geq \sum_{\tau=t_1+1}^{t_2}{(d_{k}(\tau) - f_{1k} )}.
$
Also for $k\in\mathcal{N}_D$, $b_k\geq 0$ and hence
\begin{align}\label{eq:bk>}
b_k \geq\max_{0\leq t_1< t_2\leq T} \left[ \sum_{\tau=t_1+1}^{t_2}\left(d_{k}(\tau) - f_{1k} \right)\right]^+.
\end{align}
Thus we get
$
h \geq \sum_{k\in\mathcal{N}_D}{b_k}
\geq h_{min}
$. If $\Pi^{\{1\}}(h)$ is feasible, then $P(h)$ is also feasible and hence $h\geq h_{min}$ is necessary for both problems to be feasible. Now we prove that it is also sufficient.
In particular, we show that for $h = h_{min}$, $\Pi^{\{1\}}(h)$ is feasible. For convenience, define
\begin{align}\label{eq:hk}
\tilde{h}_k := \operatorname*{max}_{0\leq t_1<t_2\leq T}{
\left[\sum_{\tau=t_1+1}^{t_2}{( d_{k}(\tau) - f_{1k})} \right]^+}, \quad k\in\mathcal{N}_D.
\end{align}
Then $h_{min} = \sum_{k \in \mathcal{N}_D}{\tilde{h}_k}$. Rearranging \eqref{eq:hk}, we get
\begin{align}\label{eq:11}
f_{1k}\geq 
\max_{0 \leq t_1<t_2 \leq T} {\left({\dfrac{\sum_{\tau=t_1+1}^{t_2} {d_k(\tau)}-\tilde{h}_k}{t_2-t_1} }\right)}  .
\end{align}
Also, by hypothesis, we have
\begin{align}\label{eq:12}
f_{1k}\geq  
\max_{t \in [T]} {\left({\dfrac{\sum_{\tau=1}^{t} {d_k(\tau)}}{t} } \right)}.
\end{align} 
Combining equations \eqref{eq:11} and \eqref{eq:12}, we have
\begin{align}
\notag
f_{1k}\geq 
\max \Big\{ &\max_{0 \leq t_1<t_2 \leq T} {\left({\dfrac{\sum_{\tau=t_1+1}^{t_2} {d_k(\tau)}-\tilde{h}_k}{t_2-t_1} }\right)},  \\
&\max_{t \in [T]} {\left({\dfrac{\sum_{\tau=1}^{t} {d_k(\tau)}}{t} } \right)}
\Big\}.\label{eq:fminrad}
\end{align}
For each $k \in \mathcal{N}_D$, consider a single generator single load system as follows. Let the demand profile be $d_k(t)$, the capacity of the transmission line be $f_{1k}$ and the total available storage budget be $\tilde{h}_k$. For this system, the right hand side in \eqref{eq:fminrad} coincides with the definition of $f_{min}$ in \eqref{eq:fmin}. From Proposition \ref{prop:SGSLf}, it follows that there is a feasible generation profile (say $g^{(k)}{(t)}$) and a storage control policy $r_k(t)$ that define a feasible flow over this single generator single load system and meet the demand. Now, for the star network, construct the generation profile $g_1(t)$
$$ g_1(t) = \sum_{k \in \mathcal{N}_D} g^{(k)}(t),$$
and operate the storage units at each node $k \in \mathcal{N}_D$ with the control policy $r_k(t)$ defined above. Also, $r_1(t) = 0$ for all $t \in [T]$. It can be checked that this defines a feasible point for $\Pi^{\{1\}}(h_{min})$.

Next, we prove that $p_*(h_{min}) = \pi^{\{1\}}_*(h_{min})$. Let $b^*_k, k\in\mathcal{N}$ be optimal storage capacities for problem $P(h_{min})$. Then the optimal storage capacities satisfy the following relations:
\begin{align*}
\sum_{k\in\mathcal{N}_D}b_k^* & \geq h_{min},\quad \text{and} \quad
b_1^* + \sum_{k\in\mathcal{N}_D} b_k^* \leq h_{min}.
\end{align*}
where the first one follows from \eqref{eq:bk>} and the second one follows from the constraint on the total available storage capacities. Rearranging the above equations, we get $b_1^*=0$ and hence $p_*(h_{min}) = \pi^{\{1\}}_*(h_{min})$. This completes the proof of part (a).


To prove part (b) of Proposition \ref{prop:hmin_starN}, we start by showing that
\begin{equation}\label{eq:infty}
p_*(\infty) = \pi^{\{1\}}_*(\infty).
\end{equation}
Assume $P(\infty)$ is feasible. For $h=\infty$, we drop the variables $b_k,k\in\mathcal{N}$, and consider the problems $P(\infty)$ and $\Pi^{\{1\}}(\infty)$ over the variables $g_1(t)$, $r_k(t), k\in\mathcal{N}$. The variables $p_{1k}(t)$ and $s_k(t)$ are defined accordingly for all $k\in\mathcal{N}$. We argue that the optimal points of $P(\infty)$ lie in a bounded set. Note that $|p_{1k}(t)| = |d_k(t) + r_k(t)| \leq f_{1k}$ and thus the control policies $r_k(t)$ are bounded for all $k\in\mathcal{N}_D$. Also, the cost function $c_1(\cdot)$ is convex and hence its sub-level sets \cite{boyd_book} are bounded. From the above arguments and the power-balance at bus 1, the optimal policy $r_1(t)$ is also bounded. Thus, the set of optimal solutions of $P(\infty)$ is a bounded set. Furthermore, this set is also closed since the objective function and the constraints are continuous functions. As in the proof of Lemma \ref{lemma:gf}, associate the function $\sum_{t\in [T]} |r_1(t)|$ with every point in the set of optimal solutions of $P(\infty)$. This is a continuous function on a compact set and hence attains a minimum. Consider the optimum of $P(\infty)$ where this minimum is attained. We prove \eqref{eq:infty} by showing that $r^*_1(t) = 0$ for all $t\in[T]$ at this optimum.


Assume to the contrary, that $r_1^*(t)$ is non-zero for some $t \in [T]$. Define
\begin{align*}
t_0 &:= \max\left\{ t\in[T] \ | \ r_1(t_0) >0 \right\},\\
t_1 &:= {\min} \left\lbrace  t\in [t_0+1,T] \ | \ r_1^*(t)<0   \right\rbrace.
\end{align*}
Also, define
$
\Delta := \min{ \left\{ r^*_1(t_0) , -r_1^*(t_1) \right\} }
$
and notice that $\Delta > 0$. 

{\emph{Case 1}}: $g_1^*(t_0) > g_1^*(t_1) +  \Delta$: Construct the modified generation and charging/ discharging profiles $\tilde{g}_1(t), \tilde{r}_1(t)$ that differ from $g_1^*(t), r_1^*(t)$ only at $t_0$ and $t_1$ as follows:
\begin{align*}
\tilde{g}_1(t_0) &= g^*_1(t_0) - \Delta g, \quad  \tilde{g}_1(t_1) = g^*_1(t_1) +  \Delta g, \\
\tilde{r}_1(t_0) &= r^*_1(t_0) - \Delta g, \quad  \ \tilde{r}_1(t_1) = r^*_1(t_1) + \Delta g,
\end{align*} 
where $\Delta g := \min\left\lbrace \Delta,g^*_1(t_0) \right\rbrace > 0$.
%
As in the proof of Lemma \ref{lemma:gf}, this is feasible for $P(\infty)$. Also, by Lemma \ref{lemma:meanVal}:
$$
 c_1\left( \tilde{g}_1(t_0)  \right) + c_1\left( \tilde{g}_1(t_1)  \right) \leq c_1\left( {g}^*_1(t_0)  \right) + c_1\left( {g}^*_1(t_1)  \right).
$$
The details are omitted for brevity.  This feasible point satisfies
\begin{align}
\label{eq:rTildeComp}
|\tilde{r}_1(t_0)| + |\tilde{r}_1(t_1)| = r_1^*(t_0) -r^*_1(t_1) -2 \Delta g  < |r_1^*(t_0)| + |r_1^*(t_1)|,
\end{align}
and hence defines an optimal point of $P(\infty)$ with a strictly lower value of the function $\sum_{t\in [T]} |r_1(t)|$. This is a contradiction.

{\emph{Case 2}}: $g_1^*(t_0) \leq g_1^*(t_1) +  \Delta$: 
As above we construct modified storage control policies $\tilde{r}_k(t)$ for all $k\in\mathcal{N}$, keeping the generation profile constant to define an optimal point of $P(\infty)$ with a lower value of $\sum_{t\in [T]} |r_1(t)|$ to derive a contradiction. 

Let the modified control policy at bus $1$ be as follows:
\begin{align*}
\tilde{r}_1(t_0) &= r^*_1(t_0) - \Delta , \quad  \tilde{r}_1(t_1) = r_1^*(t_1) + \Delta.
\end{align*} 
Instead, we distribute this to storage devices at $k\in\mathcal{N}_D$, as follows:
\begin{align*}
\tilde{r}_k(t_0) &= r^*_k(t_0) + \psi_k, \quad \tilde{r}_k(t_1) = r_k^*(t_1) -  \psi_k,\quad k\in\mathcal{N}_D,
\end{align*} 
for some $\psi_k\geq 0, k\in\mathcal{N}_D$ and $\sum_{k\in\mathcal{N}_D}{\psi_k} = \Delta$. To ensure feasibility of the modified profiles it suffices to check that the line flow constraints are satisfied at $t_0$ and $t_1$. In other words, we show that there exists $\psi_k, k\in\mathcal{N}_D$ such that for all $k\in\mathcal{N}_D$,
\begin{align*}
\psi_k \geq 0,  p_{1k}^*(t_0) + \psi_k &\leq f_{1k},  p_{1k}^*(t_1) - \psi_k\geq -f_{1k},
  \sum_{k\in\mathcal{N}_D}{\psi_k} = \Delta.
\end{align*}
Equivalently, we prove that $$
\sum_{k\in\mathcal{N}_D}{ \min\left\{ f_{1k}-p_{1k}^*(t_0) , f_{1k}+p_{1k}^*(t_1) \right\} } \geq \Delta.
$$
Recall that $p^*_{1k}(t_0)$ and $p^*_{1k}(t_1)$ are feasible for $P(\infty)$. Thus $p^*_{1k}(t_0)\leq f_{1k}$  and $p^*_{1k}(t_1)\geq -f_{1k}$. Also, $g^*_1(t) - r^*_1(t) = \sum_{k\in\mathcal{N}_D}{p^*_{1k}(t)}$ at $t=t_0$ and $t=t_1$. Thus, we have
\begin{align*}
&\sum_{k\in\mathcal{N}_D}{ \min\left\{ f_{1k}-p_{1k}^*(t_0) , f_{1k}+p_{1k}^*(t_1) \right\} } \\
&\qquad\quad\geq 
\sum_{k\in\mathcal{N}_D}{ \left( p_{1k}^*(t_1)-p_{1k}^*(t_0) \right) } \\
&\qquad\quad= \underbrace{g_1^*(t_1) - g_1^*(t_0)}_{\geq -\Delta} - \underbrace{r_1^*(t_1)}_{\leq -\Delta}   + \underbrace{r_1^*(t_0)}_{\geq \Delta}  \geq \Delta,
\end{align*}
where the last inequality follows from the hypothesis $g_1^*(t_0) \leq g_1^*(t_1) +  \Delta$. The modified profiles satisfy 
$|\tilde{r}_1(t_0)| + |\tilde{r}_1(t_1)|  <  |r_1^*(t_0)| + |r_1^*(t_1)|$ as in \eqref{eq:rTildeComp}.
As argued above this is a contradiction and hence \eqref{eq:infty} holds.

For $P(\infty)$, $s_k^*(t),k\in\mathcal{N}, t\in[T]$ is finite. Define
$
h_{sat} := \sum_{k\in\mathcal{N}_D}{ \max_{t\in[T]}s_k^*(t) }.
$
Then, note that $\left( {g}^*_1(t),r^*_k(t), t\in[T] \ k\in \mathcal{N}\right)$ are also feasible for $\Pi^{\{1\}}(h)$ and $P(h)$ for all $h\geq h_{sat}$. This
completes the proof.

\subsection{Proof of Proposition \ref{prop:proportionalStar}}
As in the statement of the proposition, assume demand profiles $d_k(t),t\in[T]$ are fixed and satisfy $d_k(t) = \phi_k d_2(t)$ with $\phi_2=1$ and $\phi_k,k=3,\ldots,n$ positive constants. Also, fix  storage budget $h$.

First, we introduce some notation. 
Denote the parametrization of $P$ with respect to the line flow capacities $f_{1k}$ as $P(\{f_{1k}\})$. For any instance of $\{f_{1k}\}$ define the problem $\Ps(f)$ as the storage placement problem in a single generator single load network with line flow constraint $f=\sum_{k=2}^{n} f_{1k}$, demand profile $d(t) = \sum_{k=2}^{n} d_k(t)$ and storage budget $h$. Denote $\ps_*(f)$ the optimal cost of $\Ps(f)$. Given an optimal solution of $P(\{f_{1k}\})$, it is straightforward to construct a feasible solution for $\Ps(f)$ with the same cost. Hence, 
\begin{align}\label{eq:from99}
\ps_*(f)\leq p_*(\{f_{1k}\}).
\end{align}

Now assume a specific realization of the line flow capacities, say $\{\hat f_{1k}\}$, such that  \eqref{eq:pro5} is satisfied. Then, it can be checked from \eqref{eq:h_minRadial} that $P(\{\hat f_{1k}\})$ is feasible for all values of the storage budget. In what follows we construct an optimal solution of $P(\{\hat f_{1k}\})$, which assigns zero storage at the generator. 
  Define
\begin{align*}
f_* =     \max_{0\leq t_1 < t_2 \leq T} \left( \frac{ \sum_{\tau=t_1+1}^{t_2} d_2(\tau)} {t_2-t_1} \right), 
\end{align*}
and  $\tilde f_{1k} := \phi_k f_*$. 
Clearly, $\{\tilde f_{1k}\}$ satisfies \eqref{eq:pro5} and $\tilde f_{1k} \leq \hat f_{1k}$. 

Let  $\hat f = \sum_{k=2}^{n}\hat f_{1k}$ and $\tilde f = \sum_{k=2}^{n}\tilde f_{1k}$. Consider  the problems $\Ps(\hat f)$ and $\Ps(\tilde f)$. It can be checked from \eqref{eq:fmin} in Proposition \ref{prop:SGSLf} that both problems are feasible. Moreover, Proposition \ref{prop:SGSLf} establishes that $\ps_*(\hat f) = \ps_*(\tilde f)$. But, $\ps_*(\hat f)\leq p_*(\{\hat f_{1k}\})$ from \eqref{eq:from99}. Combining,  
\begin{align}\label{eq:rel11}
\ps_*(\tilde f) \leq p_*(\{\hat f_{1k}\}).
\end{align} 

We start with an optimal solution of $\Ps(\tilde f)$ and construct a feasible solution with the same generation cost for $P(\{\hat f_{1k}\})$. Then, we use \eqref{eq:rel11} to deduce optimality of this solution. Our construction is such that zero storage is assigned at the generator, which establishes the desired result.

 From Theorem \ref{thm:main}, there exists an optimal solution to $\Ps(\tilde f)$ which assigns zero storage at the generator; let $g^*(t),t\in[T]$ and $r^*(t),t\in[T]$ be the optimal generation profile and optimal storage control policy of the load bus, respectively. From the feasibility constraints,
\begin{align}
g^*(t) &= \sum_{k}d_k(t) + r^*(t) \leq \tilde f = \sum_{k}\tilde f_{1k} \notag \\
&\Rightarrow\left(\sum_{k}{\phi_k}\right)d_2(t) + r^*(t) \leq \left(\sum_{k}\phi_{k}\right)f_*\notag\\
&\Rightarrow d_k(t) +  \frac{\phi_k}{\sum_{k}{\phi_k}}r^*(t) \leq \tilde f_{1k}
\end{align}
Using this, it can be checked that $g_1^*(t):=g^*(t)$, $r_k^*(t) :=  r^*(t)\phi_k/(\sum_{k}{\phi_k})$ and $b^*_k := h\phi_k/(\sum_{k}{\phi_k})$, defines a feasible solution for $P(\{\tilde f_{1k}\})$. Since $\tilde f_{1k} \leq \hat f_{1k}$ the solution is also feasible for $P(\{\hat f_{1k}\})$. In view of \eqref{eq:rel11} this completes the proof.

\end{document}